   \documentclass[11pt]{amsart}
   \usepackage{amsxtra}
   \usepackage[innercaption]{sidecap}
   \usepackage{graphicx}
   \graphicspath{{./images/}}

\theoremstyle{plain}

\def\CC {{\mathbb C}}
\def\RR {{\mathbb R}}
\def\NN {{\mathbb N}}
\def\ZZ {{\mathbb Z}}
\def\PP {{\mathbb P}}

\def\be {\begin{eqnarray}}
\def\ben {\begin{eqnarray*}}
\def\ee {\end{eqnarray}}
\def\een {\end{eqnarray*}}

\def\AAA{\kern-0.3em}
\def\AA{\kern-0.18em}
\def\AC{\kern-0.14em}
\def\AB{\kern-0.22em}

\newcommand \nc {\newcommand}

\newtheorem{theorem}{Theorem}[section]
\newtheorem{lemma}[theorem]{Lemma}
\newtheorem{proposition}[theorem]{Proposition}
\newtheorem{corollary}[theorem]{Corollary}
\newtheorem{definition}[theorem]{Definition}
\newtheorem{example}[theorem]{Example}
\newtheorem{remark}[theorem]{Remark}
\nc \bth[1] { \begin{theorem}\label{t#1} } \nc \ble[1] {
\begin{lemma}\label{l#1} } \nc \bpr[1] {
\begin{proposition}\label{p#1} } \nc \bco[1] {
\begin{corollary}\label{c#1} } \nc \bde[1] {
\begin{definition}\label{d#1}\rm } \nc \bex[1] {
\begin{example}\label{e#1}\rm } \nc \bre[1] {
\begin{remark}\label{r#1}\rm } \nc \bcon[1] {
\medskip\noindent{\it{Conjecture #1}} } \nc \bqu[1]  {
\medskip\noindent{\it{Question #1}} }
\nc {\ethe} { \end{theorem} }
 \nc {\ele} { \end{lemma} } \nc {\epr}
{ \end{proposition} } \nc {\eco} { \end{corollary} } \nc {\ede} {
\end{definition} } \nc {\eex} { \end{example} } \nc {\ere} {
\end{remark} } \nc {\econ} {\smallskip} \nc {\equ} {\smallskip}
 \nc \thref[1]{Theorem \ref{t#1}}
\nc \leref[1]{Lemma \ref{l#1}} \nc \prref[1]{Proposition
\ref{p#1}} \nc \coref[1]{Corollary \ref{c#1}} \nc
\deref[1]{Definition \ref{d#1}} \nc \exref[1]{Example \ref{e#1}}
\nc \reref[1]{Remark \ref{r#1}}
\def \B {{\mathcal B}}

\def \L {{\mathcal L}}

\def \diag { {\mathrm{diag}} }
\def \res { {\mathrm{Res}} }

 \def\AA  {\kern-0.1em}
 \def\BB  {\kern+0.1em}
 \def\BBB {\kern+0.15em}
 \def\K   {\kern+0.05em}
 \def\MK  {\kern-0.07em}
 \def\MKK {\kern-0.04em}

\begin{document}

\vspace{0.5cm}

\title[ Unfolding of the DCHE ]
{ The reducible double confluent Heun equation and a general symmetric  unfolding of  the origin }

\author[Tsvetana  Stoyanova]{Tsvetana  Stoyanova}

\date{02.03.2023}

 \maketitle

\begin{center}
{Department of Mathematics and Informatics,
Sofia University,\\ 5 J. Bourchier Blvd., Sofia 1164, Bulgaria, 
cveti@fmi.uni-sofia.bg}
\end{center}

\vspace{0.5cm}

{\bf Abstract.}
 The reducible double confluent Heun equation (DCHE) is the only DCHE whose general symmetric unfolding 
 leads to a Fuchsian equation. Contrary to the general Heun equation the unfolded Fuchsian equation has 5 singular points:
 $x_L=-\sqrt{\varepsilon}, x_R=\sqrt{\varepsilon}, x_{LL}=-1/\sqrt{\varepsilon}, x_{RR}=1/\sqrt{\varepsilon}$ and
 $x_{\infty}=\infty$. We prove that the monodromy matrix around the regular resonant singularity at the origin is realizable
 as a limit of the product of the monodromy matrices around resonant singularities $x_L$ and $x_R$ when 
 $\sqrt{\varepsilon} \to 0$ while the Stokes matrix at the irregular singularity at the origin is a limit
 of the part of the monodromy matrix around the resonant singularity $x_L$. 
 This geometrical difference between the unfolding of the two different kinds of singularities at the origin is attended with
 an analytic difference between the coefficients in the logarithmic terms including in the solution of the unfolded equation.
 While the coefficients related to the unfolding of a regular singularity have infinite limits when 
 $\sqrt{\varepsilon} \to 0$ this one related to the unfolding of an irregular singularity has a finite limit.
 We also show that the reducible
 DCHE possesses a holomorphic solution in the whole $\CC^*$ if and only if the parameters of the equation are 
 connected by a Bessel function of first kind and order depending on the non-zero characteristic exponent at the origin.



\headsep 10mm \oddsidemargin 0in \evensidemargin 0in

\section{Introduction}

        The double confluent Heun equation  (DCHE) is a second order linear ordinary differential equation having
        two irregular singular points of Poincar\'{e} rank 1 over $\CC\PP^1$. If we fix them at $x=0$ and $x=\infty$
        the standard form of DCHE writes 
         \be\label{eq}
          w'' + \left[\frac{\alpha}{x} + \frac{\beta}{x^2} + \gamma\right]\,w' + 
          \frac{\delta\,x-q}{x^2}\,w=0\,,
         \ee
         where $\alpha, \beta, \gamma, \delta$ and $q$ are arbitrary complex parameters.
        The DCHE belongs to the list of confluent Heun's equations. They were introduced and firstly studied by Decarreau
        et al. in 1978 \cite{D1, D2}.  All of them are obtained  by different confluence procedures from the general Heun equation (GHE) 
         $$\,
         w'' + \left[\frac{\alpha}{x-1} + \frac{\beta}{x} + \frac{\gamma}{x-a}\right]\,w'
         + \frac{\delta\,\epsilon\,x-q}{x\,(x-1)\,(x-a)}\,w=0\,,
         \,$$
         which is a second order Fuchsian equation with 4 singular points.
       The DCHE is obtained by a coalescence of the regular singularities
       $x=a,\,x=\infty$ and $x=0, x=1$ of the  GHE. The first confluent procedure leads to the irregular singularity
       at $x=\infty$ while the second leads to the irregular singularity $x=0$ of the DCHE (see \cite{SL}).
       The double confluent Heun equation finds many applications in superconductivity \cite{B-T},
       statistical mechanics \cite{GM}, gravity \cite{S-S}.  
        
        In this paper we apply a reverse procedure that is different from an anti-confluent procedure.
        We start with the  double confluent Heun  equation.  
		   By introducing a small complex parameter $\varepsilon$ we unfold the equation \eqref{eq}  to
               the second order  equation
                \ben
                w'' &+&
                \Big[\frac{\alpha}{2}\left(\frac{1}{x-\sqrt{\varepsilon}}
                + \frac{1}{x+\sqrt{\varepsilon}}\right) +
                \frac{\beta}{2 \sqrt{\varepsilon}} 
                \left(\frac{1}{x-\sqrt{\varepsilon}} - \frac{1}{x+\sqrt{\varepsilon}}\right) \\[0.25ex]
                &-&
                \frac{\gamma}{2 \sqrt{\varepsilon}} 
                \left(\frac{1}{\frac{1}{\sqrt{\varepsilon}}-x} +
                \frac{1}{\frac{1}{\sqrt{\varepsilon}}+x}\right)\Big ]\,w' +
                \Big[\frac{\delta}{2}
             \left(\frac{1}{x-\sqrt{\varepsilon}} + \frac{1}{x+\sqrt{\varepsilon}}\right)\\[0.2ex]
               &-&
               \frac{q}{2 \sqrt{\varepsilon}}
            \left(\frac{1}{x-\sqrt{\varepsilon}}- \frac{1}{x+\sqrt{\varepsilon}}\right) 
                 \Big]\,w=0\,.
                \een
                We call such an unfolding {\it a general symmetric unfolding}. Contrary to the
                anti-confluent procedure the general symmetric unfolding of the DCHE does not lead to a Fuchsian equation in general.
                The unfolded equation has 5 singular points. The points $x=\sqrt{\varepsilon}, x=-\sqrt{\varepsilon},
                x=1/\sqrt{\varepsilon}$ and $x=-\sqrt{\varepsilon}$ are regular singularities. When $\delta$ and $q$ are together different
                from zero the point $x=\infty$ is an irregular singularity for the unfolded equation. It becomes a regular singularity
                  if and only if $\delta=q=0$, i. e. when the DCHE is a reducible equation. 
                  This fact is the main motivation for giving our attention to  the unfolding of the reducible DCHE 
                 \be\label{eq1}
                 w'' + \left[\frac{\alpha}{x} - \frac{\beta}{x^2} - \gamma\right]\,w'=0\,,
                 \ee
                 which is obtained from the equation \eqref{eq} with $\delta=q=0$ after the transformation $x \to -x$. 
               Without  loss of generality (after a rotation of $x$) throughout this paper we assume  that $\beta$  in \eqref{eq1} is a real non-negative parameter.
                 The corresponding 
                  unfolded Fuchsian equation writes     
		   	     \be\label{eq2}
		   	      w'' &+&
		   	       \Big[\frac{\alpha}{2}\left(\frac{1}{x-\sqrt{\varepsilon}}
		   	      + \frac{1}{x+\sqrt{\varepsilon}}\right) -
		   	      \frac{\beta}{2 \sqrt{\varepsilon}} 
		   	      \left(\frac{1}{x-\sqrt{\varepsilon}} - \frac{1}{x+\sqrt{\varepsilon}}\right) \\[0.25ex]
		   	          &-&
		   	      \frac{\gamma}{2 \sqrt{\varepsilon}} 
		   	      \left(\frac{1}{\frac{1}{\sqrt{\varepsilon}}-x} +
		   	        \frac{1}{\frac{1}{\sqrt{\varepsilon}}+x}\right)\Big ]\,w'=0\,.\nonumber
		   	        \ee
		   	    We denote the singular points of the unfolded equation  by 
		   	  $x_L=-\sqrt{\varepsilon}, x_R=\sqrt{\varepsilon}, x_{LL}=-1/\sqrt{\varepsilon}, x_{RR}=1/\sqrt{\varepsilon}$
		   	  and $x_{\infty}=\infty$.
		   	  Obviously the singular points $x_L$ and $x_R$ are obtained by the unfolding of $x=0$ while 
              $x_{LL}$ and $x_{RR}$ are obtained by the unfolding of $x=\infty$ of \eqref{eq1}. It is expected that $x_{\infty}$ is
              also a result of the unfolding of $x=\infty$. In this paper comparing the analytic invariants of both equations
              we confirm this conjecture. More precisely, we will show that the analytic invariants of the DCHE around the origin 
              are realizable as a limit when $\sqrt{\varepsilon} \to 0$ of the analytic invariants of the unfolded equation only
              around resonant singularities $x_L$ and $x_R$. This phenomenon implies that the monodromy around $x_{LL}, x_{RR}$
              and $x_{\infty}$ is responsible for the unfolding of the analytic invariants around the singularity $x=\infty$
              of the DCHE. The study of the nature of the unfolding of $x=\infty$ is left to another project.
             Similar kind of problems related to the unfolding and confluence of singularities of the differential equations 
              have been studied in the works of Bolibrukh \cite{B}, Glutsyuk \cite{AG, AG1, AG2}, Hurtubise, Lambert and Rousseau
              \cite{HLR, CL-CR, CL-CR1, CL-CR2}, Klime\v{s} \cite{Kl-1, Kl-2}, Ramis \cite{R}, Stoyanova \cite{St1, St2}, Zhang {Z}.
             In the works of Buchstaber and Glutsyuk \cite{B-G}, El-Jaick and Figueiredo \cite{EJ-F}, Roseau \cite{Ro}, Tertychniy \cite{T} 
             have been studied solutions 
             space and Stokes phenomenon of the families
             of double confluent Heun equations.
                      
             The kind of singularity at the origin depends on the parameter $\beta$. When $\beta=0$ the origin is a regular
             singular point and the DCHE \eqref{eq1} degenerates into a Bessel type of equation. We introduce the notion
             of unfolded monodrmy  (see \deref{unfolded-monodromy}) as an analog of the unfolded Stokes matrix introduced by
             Lambert and Rousseay in \cite{CL-CR}. The unfolded monodromy measures geometrically the transformation
             of the monodromy around the regular singularity at the origin after a general symmetric unfolding. 
             The reducibility allows us to prove in Section 3.3, \thref{an-mon}
              that when $\beta=0$
             the monodromy around the origin of the equation \eqref{eq1} is realizable as a limit when $\sqrt{\varepsilon} \to 0$ of the 
             unfolded monodromy  which depends analytically on $\sqrt{\varepsilon}$.
             The main result in Section 3 states
             that the  monodromy matrix around the resonant singularity at the origin is realizable as a limit of
              product of the local monodromy
             matrices of the unfolded equation around resonant singular points $x_L$ and $x_R$ when $\sqrt{\varepsilon} \to 0$
             (see \prref{main1}).
             In Section 3.2, \leref{q-0} we demonstrate by a direct computation that the coefficients in the logarithmic terms
             of the solution of the unfolded equation have  limits when $\sqrt{\varepsilon} \to 0$ and both of the  limits are equal to
             $\infty$ whose sign depends on the parameter $\alpha$. It turns out that the sum of these coefficients has a finite limit
             when $\sqrt{\varepsilon} \to 0$ which is equal to the monodromy around the origin of the solution of the DCHE
             (see \coref{q}). \leref{q-0} together with \leref{d} in Section 4.1 fix the main
             difference between the unfolding of a regular singularity and an irregular singularity. In \leref{d} we show 
             explicitly that when the origin is an irregular singularity the coefficient in the logarithmic term of the
             solution of the unfolded equation has a finite limit when $\sqrt{\varepsilon} \to 0$. Moreover, this limit
             multiplied by $2\,\pi\,i$ is equal to the corresponding Stokes multiplier.
             In \cite{St2} we have shown by a direct computation that when $\alpha=2, \beta\neq 0$ the Stokes matrices at $x=0$ and $x=\infty$
              of the reducible double confluent Heun equation \eqref{eq1} are realizable as a limit of the part of the monodromy matrices
              around a resonant singularity of the general reducible Heun equation \eqref{eq2}.       
              In Section 4.1 based on the recent works of Lambert, Rousseau, Hurtubise and Klime\v{s} \cite{HLR, Kl-2, CL-CR, CL-CR1}
              we extend the result in \cite{St2} to an arbitrary reducible DCHE \eqref{eq1} without studying this equation.
              In fact this theoretical result allows us to derive the Stokes multiplier at the origin from the unfolded equation.
              In Section 4.2 we build explicit fundamental matrix solution at the origin with respect to which the
              Stokes multiplier is equal to that one obtained in Section 4.1. It turns out that the reducible DCHE \eqref{eq1}
              admits a solution which is holomorphic in whole $\CC^*$ if and only if 
              the parameters $\alpha, \beta$ and
              $\gamma$  satisfy either the relation
               \be\label{R1}
                \sum_{k=0}^{\infty} \frac{(-1)^k\,\beta^k\,\gamma^k}{k!\,\Gamma(2-\alpha+k)}=0,\quad
                \alpha\notin \NN\,,
               \ee
               or 
               \be\label{R2}
                \sum_{k=0}^{\infty} \frac{(-1)^k\,\beta^k\,\gamma^k}{k!\,\Gamma(\alpha+k)}=0,\quad
                \alpha\in \NN\,,
                \ee
             where $\Gamma(z)$ is the Euler Gamma function.
             The relations \eqref{R1} and \eqref{R2} associate the parameters $\alpha, \beta$ and $\gamma$ with the Bessel function
              $$\,
            J_{\alpha}(x)=
            \left(\frac{x}{2}\right)^{\alpha}\,
            \sum_{n=0}^{\infty}
            \frac{(-1)^n\,(x/2)^{2 n}}{n!\,\Gamma(n+\alpha+1)}
            \,$$
            of the first kind of order $\alpha$.
           
              This paper is organized as follows. In Section 2 we introduce the fundamental matrix solutions with respect to which we
              will compare the analytic invariants of both equations. We also determine the conditions on the parameters under which
              the solution of the unfolded equation can contain logarithmic terms near the singular points $x_L$ and $x_R$. 
              In Section 3 we study the unfolding of the regular singularity at the
              origin and the corresponding monodromy. The main result of Section 3 is \prref{main1} which  states 
              that when both of the singular points $x_L$ and $x_R$ are resonant singularity the monodromy matrix around the origin
              of the DCHE is realizable as a limit of the product $M_R(\varepsilon)\,M_L(\varepsilon)$ of the monodromy matrices
              around $x_R$ and $x_L$ when $\sqrt{\varepsilon} \to 0$.
              In Section 4 we deal with the unfolding of the irregular singularity
              at the origin and the corresponding Stokes phenomenon. The main result of Section 4.1 is \thref{limit}
              which states that the Stokes matrix $St_{\pi}$ at the origin of the DCHE is realizable as a limit of the part
              of the monodromy matrix around resonant singularity $x_L$ of the unfolded equation when $\sqrt{\varepsilon} \to 0$.
              The main result of Section 4.2 is \thref{0-1} which provides an actual fundamental matrix solution at
              the origin of the DCHE. The paper contains also an Appendix where we
            confirm \coref{q} by a direct computation for lower values of the parameter $\alpha$. 
              
             Since this paper appears as an extension of \cite{St2} we use without any effort some definitions and facts
             from \cite{St2}.

%%%%%%%%%%%%%%%%%%%%%%%%%%%%%%%%%%%%%%%%%
% solutions
%%%%%%%%%%%%%%%%%%%%%%%%%%%%%%%%%%%%%
         \section{Global solutions and logarithms, singular direction}

       \bth{t1}
         The equation \eqref{eq1} possesses a fundamental set of solutions $\{w_1(x, 0), w_2(x, 0)\}$ of the form
         \be\label{fss}
           w_1(x, 0)=1,\qquad
           w_2(x, 0)=\int_{\Gamma(x, 0)}
           z^{-\alpha}\,e^{-\frac{\beta}{z}}\,e^{\gamma\,z}\, d z\,.
         \ee
         The path of integration $\Gamma(x, 0)$ is  taken in such a way that the function
         $w_2(x, 0)$ is a solution of equation \eqref{eq1}.
       \ethe
       We have a similar result for the equation \eqref{eq2}.
        
        \bth{t1-p}
         The equation \eqref{eq2} possesses a fundamental set of solution $\{w_1(x, \varepsilon), w_2(x, \varepsilon)\}$
         of the form
          \be\label{fss-p}
            & &
          w_1(x, \varepsilon)=1,\\[0.1ex]
            & &
            w_2(x, \varepsilon)=\int_{\Gamma(x, \varepsilon)}
           (z-\sqrt{\varepsilon})^{\frac{\beta}{2 \sqrt{\varepsilon}}-\frac{\alpha}{2}}\,
           (z+\sqrt{\varepsilon})^{-\frac{\beta}{2 \sqrt{\varepsilon}}-\frac{\alpha}{2}}\,
            \left(\frac{\frac{1}{\sqrt{\varepsilon}}+z}{\frac{1}{\sqrt{\varepsilon}}-z}\right)^{\frac{\gamma}{2 \sqrt{\varepsilon}}}\,
            d z\,,\nonumber
            \ee
            which depends analytically on $\sqrt{\varepsilon}$.
            The path of integration $\Gamma(x, \varepsilon)$ such that $\Gamma(x, \varepsilon) \rightarrow \Gamma(x, 0)$ when
            $\sqrt{\varepsilon} \rightarrow 0$ is a path  with the same base point $x$ as the path $\Gamma(x, 0)$ 
            from \thref{t1} and taken in such a way
            that the function $w_2(x, \varepsilon)$ is a solution of the equation \eqref{eq2}.
       \ethe
        The paths $\Gamma(x, 0)$ and $\Gamma(x, \varepsilon)$ will be determined more precisely below.

       As a direct consequence of \thref{t1} and \thref{t1-p} we construct  fundamental matrices  of equations \eqref{eq1} and
       \eqref{eq2}.
       
        \bco{t2}
          The equations \eqref{eq1}  and \eqref{eq2} possess a fundamental matrix solution $\Phi(x, \cdot)$ in the form
           \be\label{fms}
             \Phi(x, \cdot)=\left(\begin{array}{cc}
             	 1   & w_2(x, \cdot)\\
             	 0   & w'_2(x, \cdot)
             	            \end{array}\right), \quad \cdot=\{0, \varepsilon\}\,,
           \ee
           where $w_2(x, \cdot), \cdot=\{0, \varepsilon\}$ is defined by \thref{t1} and \thref{t1-p}, respectively.
        \eco

        Let us determine when the solution $w_2(x, \varepsilon)$ of the unfolded equation can contain logarithmic terms near
        the singular points $x_j, j=L, R$. Recall that from the local theory of the Fuchsian singularity such a singular point
        is called a resonant singularity.  When $\beta=0$ the points $x_L$ and $x_R$ are together
        either non-resonant or resonant singularities for the unfolded equation. In particular, they both are resonant
        singularities if and only if $\alpha\in 2 \NN$.  In the next  section we consider the equations \eqref{eq1} and
        \eqref{eq2} under the restriction
        \be\label{res1}
        \beta=0,\qquad \alpha\in 2 \NN\,.
        \ee   
        Note that under the restriction \eqref{res1} the origin is a resonant regular singularity too. 
          Using the rotation $x \to x\,e^{i \delta}$ where $\delta=\arg (\sqrt{\varepsilon})$ we always can fix $\sqrt{\varepsilon}$
          to be a real and positive. Due to this property when $\beta=0$ we choose  the path $\Gamma(x, 0)$ in \eqref{fss} to be 
          a path from
          $1$ to $x$ approaching 1 in the direction $\RR^{+}$.  The path $\Gamma(x, \varepsilon)$  is a path taken in the same direction
          $\RR^{+}$
          from $1+\sqrt{\varepsilon}$ to the same base poin $x$.

        When $\beta > 0$  we choose the path $\Gamma(x, 0)$ in \eqref{fss} to be a path from $0$ to $x$ approaching $0$ in the direction 
        $\RR^{+}$. Then  
        the corresponding unfolded path $\Gamma(x, \varepsilon)$ is  a path  taken
        in the same direction $\RR^{+}$ from $\sqrt{\varepsilon}$ to the same base point $x$. This choice of the path
        $\Gamma(x, \varepsilon)$ implies that $\varepsilon$ is a real  positive parameter of unfolding and that $x_L$ will be the resonant singularity. 
        In particular in Section 4 we consider the unfolded equation under the restriction
        \be\label{res2}
         \frac{\beta}{2 \sqrt{\varepsilon}} + \frac{\alpha}{2}\in\NN,\quad
         \frac{-\beta}{2 \sqrt{\varepsilon}} + \frac{\alpha}{2}\notin\NN\,.
         \ee
      We denote by $\Phi_0(x, 0)$ and $\Phi_0(x, \varepsilon)$ the fundamental matrix solutions from \eqref{fms} 
      corresponding to the so chosen paths $\Gamma(x, 0)$ and $\Gamma(x, \varepsilon)$.
         
      From Definition 6.15 in \cite{St2} it follows that $\theta=\arg (0-\beta)=\arg(-\beta)=\pi$ is the only possible
      singular direction at the origin of the DCHE.

        \section{The unfolding of the monodromy around the origin}
        
         In this section we deal with the equations \eqref{eq1} and \eqref{eq2} when the parameters $\alpha$ and $\beta$ satisfy the
         condition \eqref{res1}.

         \subsection{The monodromy around the origin of the DCHE }

            Since the origin is a regular point for the DCHE  its unfolding causes an unfolding of the monodromy around it.
            To compute this monodromy we rewrite the fundamental matrix $\Phi_0(x, 0)$ in an appropriate form.
           Directly from \eqref{fss} and \eqref{fms} we have

        	 \bth{local-0}
        	 Assume that the condition \eqref{res1} holds. Then  the fundamental matrix solution $\Phi_0(x, 0)$ of the equation
        	 \eqref{eq1} is represented in a neighborhood of the origin as
        	 \be\label{sol-0}
        	 \Phi_0(x, 0)=\exp(G x)\,\,H(x)\,x^{\Lambda}\,x^J\,,
        	 \ee
        	 where
        	 $$\,
        	 G=\diag (0, \gamma),\quad 
        	 \Lambda=\diag(0, -\alpha)\,.
        	 \,$$
        	 The matrix $H(x)$ is defined as
        	 $$\,
        	 H(x)=\left(\begin{array}{cc}
        	 1   & x\,\varphi(x)\\
        	 0   & 1
        	 \end{array}\right)\,,
        	 \,$$            	 
        	 where $\varphi(x)$ is a holomorphic function in a neighborhood of the origin.
        	 The matrix $J$ is given by
        	 $$\,
        	 J=\left(\begin{array}{cc}
        	 0    &\lambda\\
        	 0    &0
        	 \end{array}\right)\,,
        	 \,$$
        	 where
        	\be\label{a}
        	 	\lambda=
        	 \frac{\gamma^{\alpha-1}}{(\alpha-1)!}\,.
        	 	\ee

                   \ethe

                   The monodromy of the fundamental matrix solution $\Phi_0(x, 0)$ around the origin is described
                   by the local monodromy matrix $M_0\in GL_2(\CC)$
                    \be\label{mm-0}
                    M_0=e^{2 \pi\,i\,\Lambda}\,e^{2 \pi\,i\,J}=e^{2 \pi\,i\,J}=
                    \left(\begin{array}{cc}
                    	1   & 2 \pi\,i\,\lambda\\
                    	0   & 1
                    \end{array}
                    \right)\,,
                    \ee
                    where $\lambda$ is itroduced by \eqref{a}.

         \subsection{The monodromy around $x_L$ and $x_R$ of the unfolded equation}

        In this section we compute the local monodromy matrices of the equation \eqref{eq2} under the restriction \eqref{res1}.

        In the next theorem we describe the local behavior of the fundamental matrix $\Phi_0(x, \varepsilon)$ near the singular points 
        $x_R$ and $x_L$ when both of them are resonant singularities.
        \bth{F-B-0}
        Assume that the condition \eqref{res1} holds.
        Then the fundamental matrix solution $\Phi_0(x, \varepsilon)$ of 
        the unfolded
        equation depends analytically on $\sqrt{\varepsilon}$ and it is represented in a neighborhood of the origin which  contains only the  singular points
        $x_L$ and $x_R$ as
        \ben
        \Phi_0(x, \varepsilon) =
        G(x, \varepsilon)\,H(x, \varepsilon)\,
        (x-x_L)^{\frac{1}{2} \Lambda}\,
        (x-x_R)^{\frac{1}{2} \Lambda}\,
        (x-x_L)^{J_L(\varepsilon)}\, (x-x_R)^{J_R(\varepsilon)}\,,    
        \een
        where 
        $$\,
        G(x, \varepsilon)=(x-x_{LL})^{-\frac{x_{LL}}{2}\,G}\, (x_{RR}-x)^{-\frac{x_{RR}}{2}\,G}\,.
        \,$$
        The matrix $H(x, \varepsilon)$ is a holomorphic matrix-function at the both singular
        points $x_L$ and $x_R$ such that $H(x_k, \varepsilon)=I_2, k=L, R$. The matrices $G$ and $\Lambda$ are introduced in
        \thref{local-0}. The matrix $J_k(\varepsilon), k=L, R$ is given by
        \be\label{J-B-0}
        J_k(\varepsilon)=\left(\begin{array}{cc}
        	0   &q_k\\
        	0   &0
        \end{array}
        \right)\,,
        \ee
        where the number $q_k$  is defined as 
        \be\label{q-B-0}
        q_k=\res \left( w_2'(x, \varepsilon),\,x=x_k\right),\quad
        k=L, R\,.
        \ee
        \ethe     
        
        \proof
        
        Let us  present the fundamental matrix $\Phi_0(x, \varepsilon)$
        in the form
        \ben
        & &
        \Phi_0(x, \varepsilon)=
        \left(\begin{array}{cc}
        	1   & 0\\
        	0                       & w_2'(x, \varepsilon)
        \end{array}
        \right)
        \left(\begin{array}{cc}
        	1   &\int_{1+\sqrt{\varepsilon}}^x w'_2(z, \varepsilon)\,d z\\
        	0   &1
        \end{array}
        \right)\\[0.2ex]
        &=&
        G(x, \varepsilon)\,\left[(x-x_L) (x-x_R)\right]^{\frac{1}{2} \Lambda}
        \left(\begin{array}{cc}
        	1   &\int_{1+\sqrt{\varepsilon}}^x w'_2(z, \varepsilon)\,d z\\
        	0   &1
        \end{array}
        \right)\,,    	  
        \een
        where
        $$\,
        G(x, \varepsilon)=(x-x_{LL})^{-\frac{x_{LL}}{2}\,G}\, (x_{RR}-x)^{-\frac{x_{RR}}{2}\,G}\,.
        \,$$

        Consider the function
        $ w_2 (x, \varepsilon)$. Since when $\alpha\in 2 \NN$ the function 
        $\frac{1}{[(z-\sqrt{\varepsilon})(z+\sqrt{\varepsilon})]^{\alpha/2}}$ is a rational function it can be splited into 
        a finite sum in $\alpha$ number simpler ratios $\frac{c_j}{(z-\sqrt{\varepsilon})^j}$ and
        $\frac{d_j}{(z+\sqrt{\varepsilon})^j},\,1 \leq j \leq \alpha/2$ where the coefficients $c_j, d_j$ are uniquely
        determined. Then the function $w_2(x, \varepsilon)$  can be written
        as
        \ben
        & &
        \int_{1+\sqrt{\varepsilon}}^x w'_2(z, \varepsilon)\,d z=
        \int_{1+\sqrt{\varepsilon}}^x \frac{P(z, \varepsilon)}{(z -\sqrt{\varepsilon})^{\frac{\alpha}{2}}}
        \left(\frac{\frac{1}{\sqrt{\varepsilon}}-z}{\frac{1}{\sqrt{\varepsilon}}+z}\right)
        ^{\frac{\gamma}{2 \sqrt{\varepsilon}}}
                d z\\[0.25ex]
        &+&
        \int_{1+\sqrt{\varepsilon}}^x \frac{Q(z, \varepsilon)}{(z +\sqrt{\varepsilon})^{\frac{\alpha}{2}}}
        \left(\frac{\frac{1}{\sqrt{\varepsilon}}+z}{\frac{1}{\sqrt{\varepsilon}}-z}\right)
        ^{\frac{\gamma}{2 \sqrt{\varepsilon}}}
                d z\\[0.25ex]
        &=&
        q_R \log (x-x_R) + q_L \log (x-x_L) + (x-x_R)^{-\frac{\alpha}{2}+1}\,h(x-x_R) + (x-x_L)^{-\frac{\alpha}{2}+1}\,g(x-x_L)\,.
        \een         
        Here $P(z, \varepsilon)$ and $Q(z, \varepsilon)$ are polynomials of degree at most $\alpha/2-1$.        
        The functions $h(x-x_R)$ and $g(x-x_L)$ are holomorphic functions at the both singular points $x_R$ and $x_L$ since the
        function
        $\left(\frac{\frac{1}{\sqrt{\varepsilon}}+z}{\frac{1}{\sqrt{\varepsilon}} -z}\right)
        ^{\frac{\gamma}{2 \sqrt{\varepsilon}}}$ is a holomorphic function at the both singular 
        points $x_j, j=L, R$. 
        Then we can present
        $\Phi_0(x, \varepsilon)$ as
        \ben
        & &
        G(x, \varepsilon)
        \left(\begin{array}{cc}
        	1  &(x-x_R)(x-x_L)^{\frac{\alpha}{2}} h(x-x_R) + (x-x_L)(x-x_R)^{\frac{\alpha}{2}} g(x-x_L)\\
        	0  &1
        \end{array}\right)\\[0.25ex]
        &\times&
        \left(\begin{array}{cc}
        	1    &q_L \log (x-x_L) + q_R \log (x-x_R)\\
        	0    &(x-x_L)^{-\frac{\alpha}{2}} (x-x_R)^{-\frac{\alpha}{2}}
        \end{array}\right)\\[0.25ex]
        &=&
        G(x, \varepsilon)
        \left(\begin{array}{cc}
        	1  &(x-x_R) (x-x_L)^{\frac{\alpha}{2}} h(x-x_R) + (x-x_L) (x-x_R)^{\frac{\alpha}{2}} g(x-x_L)\\
        	0  &1
        \end{array}\right)\\[0.25ex]
        &\times&
        [(x-x_L) (x-x_R)]^{\frac{1}{2} \Lambda}     			
        \left(\begin{array}{cc}
        	1    &q_L \log (x-x_L) + q_R \log (x-x_R)\\
        	0    &1
        \end{array}\right)\\[0.2ex]
        &=& 	 
        G(x, \varepsilon)
        \left(\begin{array}{cc}
        	1  &(x-x_R) (x-x_L)^{\frac{\alpha}{2}} h(x-x_R) + (x-x_L) (x-x_R)^{\frac{\alpha}{2}} g(x-x_L)\\
        	0  &1
        \end{array}\right)\\[0.25ex]
        &\times&
        [(x-x_L) (x-x_R)]^{\frac{1}{2} \Lambda}     			
        (x-x_L)^{J_L(\varepsilon)} (x-x_R)^{J_R(\varepsilon)}=\Phi_0(x, \varepsilon)\,.
        \een
        This ends the proof.
        \qed

        Consider the DCHE \eqref{eq1} and its fundamental matrix solution at the origin in the
        punctured disk $D_R$ around the origin with a finite small radius $R$
        $$\,
        D_R:=\{x\in\CC\,|\, 0 < |x| < R\}\,.
        \,$$ 
        The radius $R$ is so chosen that the points $x_L$ and $x_R$ belong to $D_R$ while the points $x_{LL}$ and $x_{RR}$
        do not belong to $D_R$. 
         Let   $x_0\in D_R \backslash \RR$. Let $\gamma_L$ and $\gamma_R$ be two closed loops, starting and ending at the
         point $x_0$. The loop $\gamma_L$ (resp. $\gamma_R$) encircles only the point $x_L$ (resp. $x_R$) in the positive sense
         as it is shown in Figure \ref{fig:image}.        
        Thanks to \thref{F-B-0}  we can fix explicitly the corresponding local monodromy matrices corresponding to the loops
        $\gamma_L$ and $\gamma_R$.
        
        \bth{M-B-0}
       The local monodromy matrices $M_k(\varepsilon),\,k=R, L$ of the perturbed equation
        with respect to the fundamental matrix $\Phi_0(x, \varepsilon)$, introduced by 
        \thref{F-B-0} are
        given by
        \be\label{M-0-B}
        M_k(\varepsilon)=e^{2 \pi\,i\,J_k(\varepsilon)}, \quad k=L, R\,.
        \ee 
        \ethe
        
        \proof
        Analytic continuation of the fundamental matrix $\Phi_0(x, \varepsilon)$ along the loop $\gamma_L$ leads to the
        new fundamental matrix $[\Phi_0(x, \varepsilon)]_{\gamma_L}$. The connection between these two fundamental matrices
        is measured geometrically by the monodromy matrix $M_L(\varepsilon)$
        $$\,
        [\Phi_0(x, \varepsilon)]_{\gamma_L}=
        \Phi_0(x, \varepsilon)\,M_L(\varepsilon)\,.
        \,$$
        On the other hand thanks to the explicit local form of the matrix $\Phi_0(x, \varepsilon)$ from \thref{F-B-0}
        we find that
        \ben
        [\Phi_0(x, \varepsilon)]_{\gamma_L}
        &=&
        G(x, \varepsilon)\, H(x, \varepsilon)\,
        (x-x_L)^{\frac{1}{2} \Lambda }\, e^{\pi\,i \Lambda}\\[0.2ex]
        &\times&
        (x-x_R)^{\frac{1}{2} \Lambda }\,
        (x-x_L)^{J_L(\varepsilon)} \,e^{2 \pi\,i\,J_L(\varepsilon)}\, (x-x_R)^{J_R(\varepsilon)}\\[0.2ex]
        &=&
        G(x, \varepsilon)\,H(x, \varepsilon)\,
        (x-x_L)^{\frac{1}{2} \Lambda}\\[0.2ex]
        &\times&
        (x-x_R)^{\frac{1}{2} \Lambda}
        (x-x_L)^{J_L(\varepsilon)} e^{2 \pi\,i\,J_L(\varepsilon)} 
        (x-x_R)^{J_R(\varepsilon)}  \,,    
        \een
        since  $e^{\pi\,i \Lambda}=I_2$ when $\alpha\in 2 \NN$.
        As a result we have 
        \ben
        [\Phi_0(x, \varepsilon)]_{\gamma_L}=
        \Phi_0(x, \varepsilon) e^{2 \pi\,i\,J_L(\varepsilon)} \,,    
        \een
        since the matrices $e^{2 \pi\,i\,J_L(\varepsilon)}$ and $(x-x_R)^{J_R(\varepsilon)}$ commute.
        In the same manner one can derive the formula for the monodromy matrix $M_R(\varepsilon)$.
        
        This ends the proof.
        \qed

        In what follows we compute the numbers $q_k, k=L, R$.
        
        \bpr{q-0}
        Assume that the condition \eqref{res1} holds.        
        Then  the number $q_R$ is given by
        \ben
        q_R=\frac{1}
        	{\left(\frac{\alpha-2}{2}\right)!}
        \sum_{k=0}^{\frac{\alpha-2}{2}}
        \left( \left(\begin{array}{c}
        	\frac{\alpha-2}{2}\\
        	k
        \end{array}\right)
        (-1)^{\frac{\alpha-2}{2}-k}
        \frac{\Gamma(\alpha-1-k)}{\Gamma(\frac{\alpha}{2})} (2 \sqrt{\varepsilon})^{-\alpha+1+k}\right)\,A_R\,,
        \een 
        where
        $$\,
        A_R=\left(\frac{\sqrt{\varepsilon}}{1+\varepsilon}\right)^k\,
        \left(\frac{1+\varepsilon}{1-\varepsilon}\right)^{\frac{\gamma}{2\sqrt{\varepsilon}}}
        \sum_{s=0}^k \left(\begin{array}{c}
        k\\
        s
        \end{array}\right)
       \frac{\Gamma(\frac{\gamma}{2 \sqrt{\varepsilon}} + s)\,\Gamma(\frac{\gamma}{2 \sqrt{\varepsilon}}+1)}
       {\Gamma(\frac{\gamma}{2 \sqrt{\varepsilon}})\,\Gamma(\frac{\gamma}{2 \sqrt{\varepsilon}}-k+s+1)}
        \left(\frac{1+\varepsilon}{1-\varepsilon}\right)^s\,.
        \,$$	
        
        Similarly,  the number $q_L$ is 
        given by
        \ben
        q_L=\frac{1}
        	{\left(\frac{\alpha-2}{2}\right)!}
        \sum_{k=0}^{\frac{\alpha-2}{2}}
        \left( \left(\begin{array}{c}
        	\frac{\alpha-2}{2}\\
        	k
        \end{array}\right)
        (-1)^{\frac{\alpha-2}{2}-k}
        \frac{\Gamma(\alpha-1-k)}{\Gamma(\frac{\alpha}{2})} (-2 \sqrt{\varepsilon})^{-\alpha+1+k}\right)\,A_L\,,
        \een 
        where
       \be\label{al}\qquad\quad
        A_L=\left(\frac{\sqrt{\varepsilon}}{1-\varepsilon}\right)^k
        \left(\frac{1-\varepsilon}{1+\varepsilon}\right)^{\frac{\gamma}{2\sqrt{\varepsilon}}}
        \sum_{s=0}^k \left(\begin{array}{c}
        k\\
        s
        \end{array}\right)
       \frac{\Gamma(\frac{\gamma}{2 \sqrt{\varepsilon}} + s)\,\Gamma(\frac{\gamma}{2 \sqrt{\varepsilon}}+1)}
       	{\Gamma(\frac{\gamma}{2 \sqrt{\varepsilon}})\,\Gamma(\frac{\gamma}{2 \sqrt{\varepsilon}}-k+s+1)}
        \left(\frac{1-\varepsilon}{1+\varepsilon}\right)^s\,.
        \ee
        \epr
        
        It turns out the numbers $q_R$ and $q_L$ given by \prref{q-0} have a limit when $\sqrt{\varepsilon}\rightarrow 0$.
        	
        \ble{q-0}
         Assume that the condition  \eqref{res1} holds. Then for each fixed $\alpha$ the numbers $q_R$ and $q_L$ computed by
         \prref{q-0} satisfy the limits
           $$\,
             \lim_{\sqrt{\varepsilon} \rightarrow 0} q_R=-
              (-1)^{\frac{\alpha}{2}}\,\infty,\qquad
                 \lim_{\sqrt{\varepsilon} \rightarrow 0} q_L=
                 (-1)^{\frac{\alpha}{2}}\,\infty
           \,$$
           when $\sqrt{\varepsilon} \rightarrow 0$.
        \ele
        
        \proof
        Applying the limit
           \be\label{l}
             \lim_{|z| \rightarrow \infty}
             \frac{\Gamma (z+a)}{\Gamma(z)\,z^a}=1
           \ee
           we find that 
            \ben
             q_R  &\to&
              -\left(\frac{1+\varepsilon}{1-\varepsilon}\right)^{\frac{\gamma}{2 \sqrt{\varepsilon}}}
             \frac{(-1)^{\frac{\alpha}{2}}}{(2 \sqrt{\varepsilon})^{\alpha-1}\,\Gamma(\frac{\alpha}{2})}
            \sum_{k=0}^{\frac{\alpha-2}{2}}
             \frac{(-1)^k\,\Gamma(\alpha-1-k)}{k!\,(\frac{\alpha-2}{2}-k)!} 
             \left(\frac{2 \sqrt{\varepsilon}\,\gamma}
                  {1-\varepsilon^2}\right)^k,\\[0.2ex]
                    q_L  &\to& 
                     \left(\frac{1-\varepsilon}{1+\varepsilon}\right)^{\frac{\gamma}{2 \sqrt{\varepsilon}}}
                    \frac{(-1)^{\frac{\alpha}{2}}}{(2 \sqrt{\varepsilon})^{\alpha-1}\,\Gamma(\frac{\alpha}{2})}
                    \sum_{k=0}^{\frac{\alpha-2}{2}}
                    \frac{\Gamma(\alpha-1-k)}{k!\,(\frac{\alpha-2}{2}-k)!} 
                    \left(\frac{2 \sqrt{\varepsilon}\,\gamma}
                    {1-\varepsilon^2}\right)^k  
            \een
            when $\sqrt{\varepsilon} \to 0$.
           Now the statement follows from the observation that $q_R$ and $q_L$ are expressed as finite sums and from the
           limits
           \be\label{ll}
             \lim_{\sqrt{\varepsilon} \rightarrow 0}
             \left(\frac{1+\varepsilon}{1-\varepsilon}\right)^{\frac{\gamma}{2 \sqrt{\varepsilon}}}=1,\quad
              \lim_{\sqrt{\varepsilon} \rightarrow 0}
              \left(\frac{1-\varepsilon}{1+\varepsilon}\right)^{\frac{\gamma}{2 \sqrt{\varepsilon}}}=1\,. 
            \ee
        \qed
        
        From \leref{q-0} it follows that the sign of the limit of the number $q_k, k=L, R$ depends on the parameter
        $\alpha$ but we always have that
          $$\,
           \lim_{\sqrt{\varepsilon} \rightarrow 0} (q_R + q_L)=\infty - \infty\,.
          \,$$

        \bre{unfolded1}
        The result of \leref{q-0} is the identification mark of the unfolding of a resonant regular point. 
        Recall that in our previous works \cite{St1, St2} all the coefficients place before the logarithmic terms in the solution
        of the unfolded equation have a finite limit  when $\sqrt{\varepsilon} \rightarrow 0$. 
        But in all previous cases these logarithmic terms measure how the Stokes matrices of the initial equation are transformed
        to the monodromy matrices of the unfolded equation. This time the logarithmic terms in the solution of the unfolded
        equation correspond to an unfolding of the monodromy matrix of the initial equation to two monodromy matrices of the 
        unfolded equation.
        \ere

	%%%%%%%%%%%%%%%%%%%%%%%%%%%%%%%%%%%%%%%%%%%%%%%%%%%%%%
	% Main results
	%%%%%%%%%%%%%%%%%%%%%%%%%%%%%%%%%%%%%%%%%%%%%%%%
	\subsection{ The unfolded monodromy around the origin  }

   In this section we connect by a radial limit $\sqrt{\varepsilon} \rightarrow 0$ the monodromy 
   matrices $M_j(\varepsilon), j=R, L$  of the unfolded equation with the monodromy matrix $M_0$ 
  of the DCHE. 
  The following proposition is a key for our study.

  \bpr{con}
  When $\sqrt{\varepsilon} \to 0$ the fundamental set of solutions $\{w_1(x, \varepsilon), w_2(x, \varepsilon)\}$
  of the unfolded equation fixed by \thref{t1-p} converges uniformly on compact sets of $D_R$ to the fundamental set of solutions
  $\{w_1(x, 0), w_2(x, 0)\}$ of the DCHE fixed by \thref{t1}.
  \epr
  
  Thanks to  \prref{con}  we have the following property of the fundamental matrices
  $\Phi_0(x, 0)$ and $\Phi(x, \varepsilon)$.
  \bco{conver}
  The fundamental matrix $\Phi_0(x, \varepsilon)$ of the unfolded equation given by \thref{F-B-0} converges uniformly on compact 
  sets of $D_R$ to the fundamental matrix $\Phi_0(x, 0)$ of the DCHE given by \thref{local-0}
  when $\sqrt{\varepsilon} \rightarrow 0$.
  \eco

   \begin{figure}[t]
   	
   	\begin{minipage}[c][1\width]
   		{0.003\textwidth}
   		\centering
   		\includegraphics[scale=0.38]{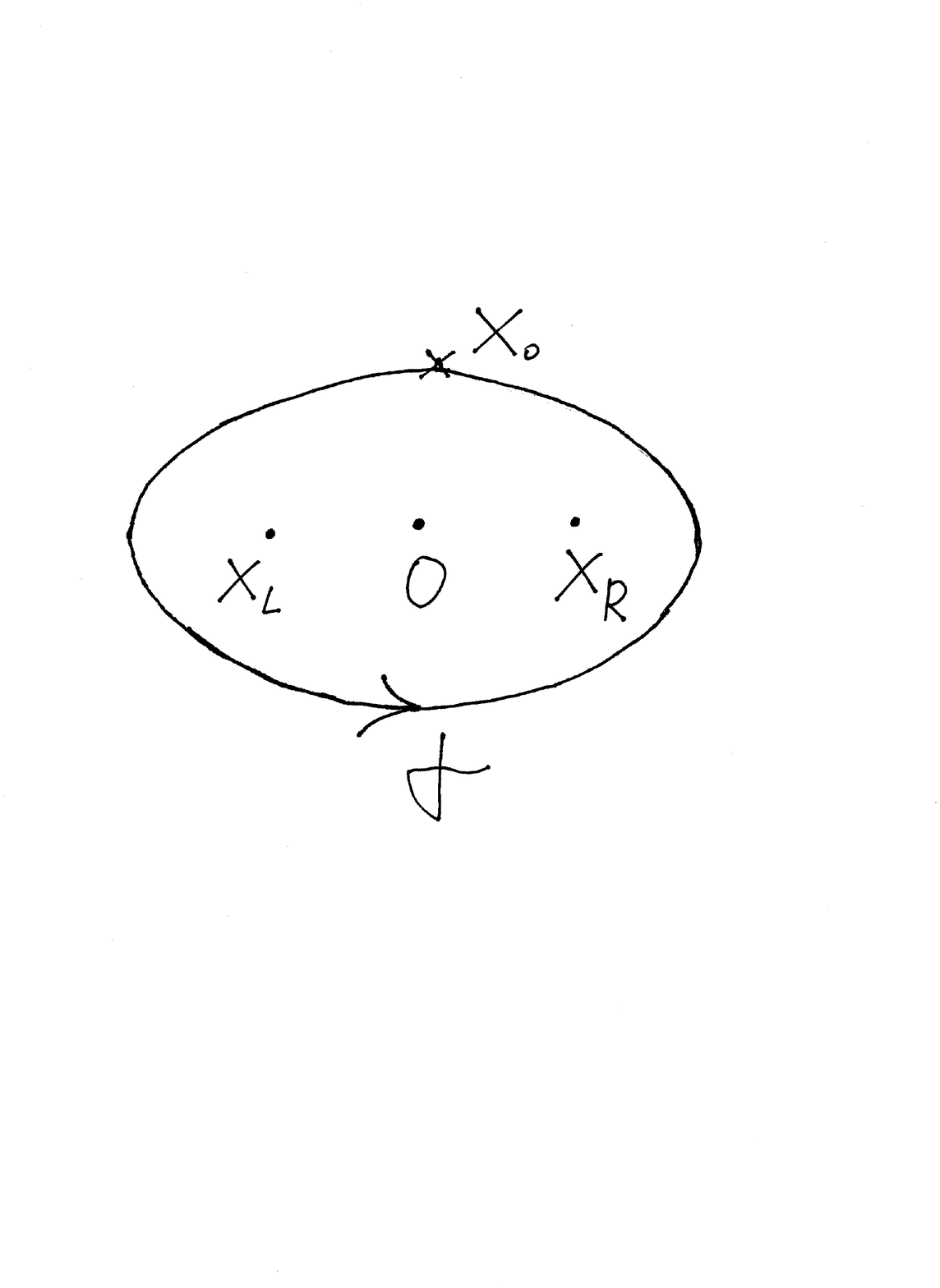}
   	\end{minipage}
   	\hfill
   	\begin{minipage}[c][1\width]
   		{0.5\textwidth}
   		\centering
   		\includegraphics[scale=0.4]{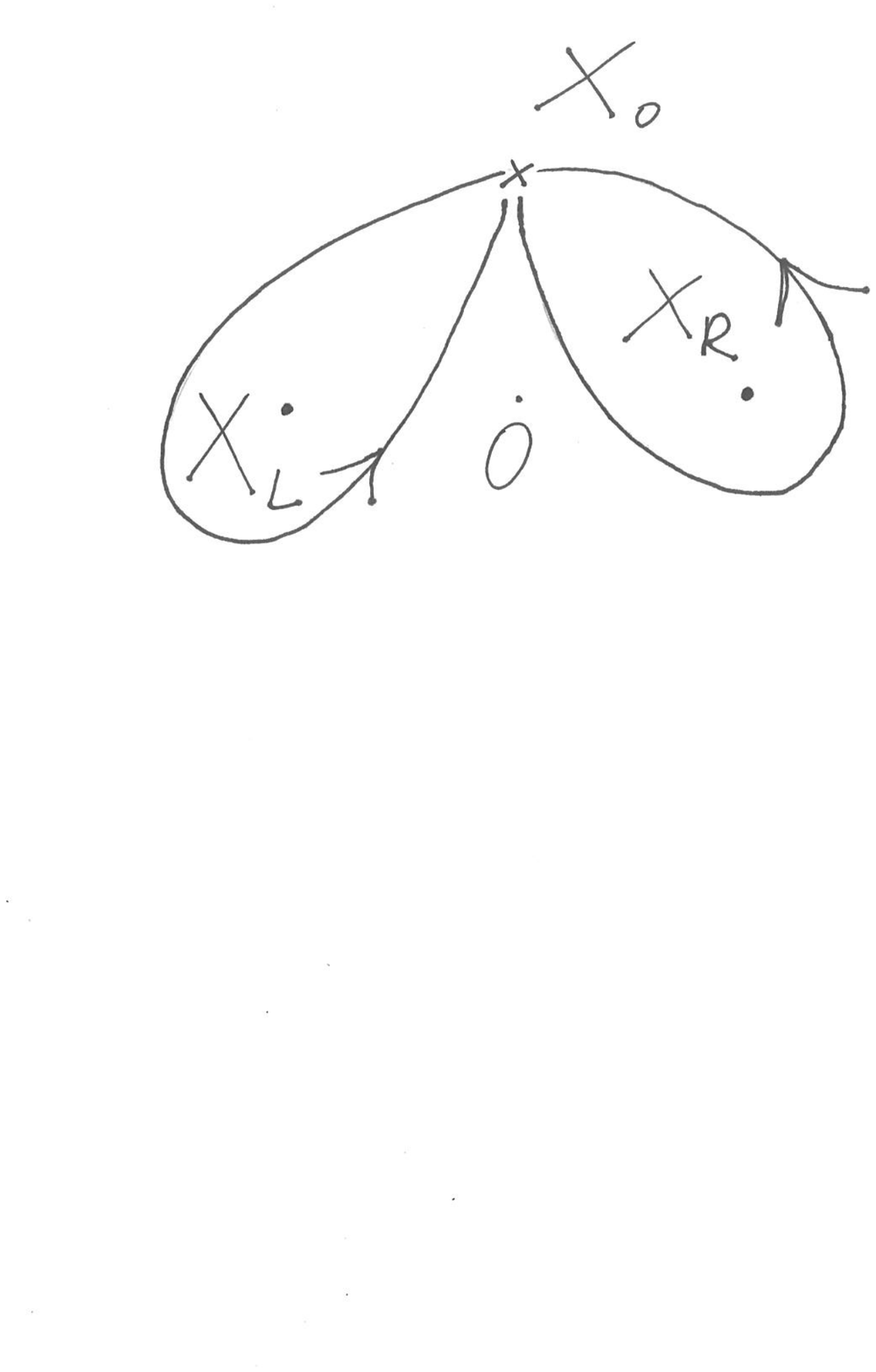}
   	\end{minipage}
   	\caption{The loop $\gamma$ in $D_R$ and the loops $\gamma_L$ and $\gamma_R$ in $D_R(\varepsilon)$.}
   	\label{fig:image}
   \end{figure}

 Let $\gamma\in D_R$ be a closed loop starting and ending at the same point 
    $x_0\in D_R \backslash \RR$ as in Section 3.2, 
    encircling the origin and the points $x_R$ and $x_L$ and oriented counter-clockwise as in  Figure \ref{fig:image}.
    The fundamental matrix solution $\Phi_0(x, 0)$ of the DCHE introduced in \thref{local-0}
    is a holomorphic multi-valued function on $D_R$. The analytic continuation of $\Phi_0(x, 0)$ along
    $\gamma$ leads to a new fundamental matrix solution $[\Phi_0(x, 0)]_{\gamma}$ of the DCHE.
    The connection between these two fundamental matrix solutions is given by the  monodromy matrix 
    $M_0=e^{2 \pi\,i\,\Lambda}\,e^{2 \pi\,i\,J}=e^{2 \pi\,i\,J}$ from
    \eqref{mm-0}
     $$\,
      \left[\Phi_0(x, 0)\right]_{\gamma}=\Phi_0(x, 0)\,M_0\,.
     \,$$
    
    Let $D_R(\varepsilon)$ be a domain in $\CC\setminus \{x_L, x_R\}$ such that  
     $D_R(\varepsilon)$ tends to the disk $D_R$ when $\sqrt{\varepsilon} \to 0$ .
    Let $\gamma(\varepsilon)\in D_R(\varepsilon)$ be a closed  loop starting and ending with the same point $x_0$ and such that 
    $\gamma(\varepsilon)=\gamma_L \circ \gamma_R$ where the loops $\gamma_L$ and $\gamma_R$ are defined in Section 3.2 and in
    Figure \ref{fig:image}.
    Then when $\sqrt{\varepsilon} \to 0$ the loop $\gamma(\varepsilon)$ tends to a closed loop that
   belongs to the homotopy class $[\gamma]$ of the loop. 
   Analytic continuation of the fundamental matrix $\Phi_0(x, \varepsilon)$ along  $\gamma(\varepsilon)$
     yields a new fundamental matrix $\left[\Phi_0(x, \varepsilon)\right]_{\gamma(\varepsilon)}$.
    The connection between these two fundamental matrices is measured  geometrically by an invertible constant
    matrix $M_0(\varepsilon)$
     \be\label{unf-mon}
       \left[\Phi_0(x, \varepsilon)\right]_{\gamma(\varepsilon)}=\Phi_0(x, \varepsilon)\,M_0(\varepsilon)\,.
     \ee
     
     \bde{unfolded-monodromy}
     We call the invertible matrix $M_0(\varepsilon)$ defined by \eqref{unf-mon} the unfolded  monodromy matrix
     around the origin.
     \ede

     	The reducibility ensures the connection by a limit $\sqrt{\varepsilon} \to 0$ between the fundamental matrix solutions
     $\Phi_0(x, 0)$ and $\Phi_0(x, \varepsilon)$, which ensures such a connection between the monodromy around the origin
     and the unfolded monodrmy
     \bth{an-mon}
      The unfolded  monodromy matrix $M_0(\varepsilon)$ around the origin depends analytically on $\sqrt{\varepsilon}$
      and converges when $\sqrt{\varepsilon} \rightarrow 0$ to the  monodromy matrix around the origin $M_0$
      defined by \eqref{mm-0}.
     \ethe
     
     \proof
      Since the fundamental matrix $\Phi_0(x, \varepsilon)$ converges uniformly on the compact sets 
      of $D_R$
      to the fundamental matrix $\Phi_0(x, 0)$, so does the fundamental matrix 
      $\left[\Phi_0(x, \varepsilon)\right]_{\gamma(\varepsilon)}$ to the fundamental matrix 
      $\left[\Phi(x, 0)\right]_{\gamma}$. 
     Then the matrix $M_0(\varepsilon)$ must converge to the  monodromy matrix $M_0$ when 
      $\sqrt{\varepsilon} \rightarrow 0$.
       \qed
     
      When the origin is a resonant singularity we find 
       that the unfolded  monodromy matrix $M_0(\varepsilon)$ is expressed in terms of the monodromy
       matrices $M_L(\varepsilon)$ and $M_R(\varepsilon)$.

    \bth{un}     
      	Assume that the condition \eqref{res1} holds.  
    Let $M_j(\varepsilon), j=L, R$ and $M_0(\varepsilon)$ be the monodromy matrices and the unfolded
     monodromy  matrix of the unfolded equation with respect to the fundamental matrix solution $\Phi_0(x, \varepsilon)$
    Then  they 
    satisfy the following relation
      \be\label{M-hat}
       M_0(\varepsilon)=M_R(\varepsilon)\,M_L(\varepsilon)=M_L(\varepsilon)\,M_R(\varepsilon)\,.
      \ee
      \ethe
   
   \proof
   The connection   $M_0(\varepsilon)=M_R(\varepsilon)\,M_L(\varepsilon)$
    follows from the definition  $\gamma_L \circ \gamma_R=\gamma(\varepsilon)$, where $\gamma_L$ and $\gamma_R$ are the loops
   from section 3.2. 
   The equality  $M_R(\varepsilon)\,M_L(\varepsilon)=M_L(\varepsilon)\,M_R(\varepsilon)$ follows from the fact that 
   under the  condition \eqref{res1}  the matrices $M_{L}(\varepsilon)$ and $M_R(\varepsilon)$ commute.
   \qed

    	As an immediate consequence we have
    	\bco{J}
    	Assume that the condition \eqref{res1} holds. 
    	Then the unfolded  monodromy  matrix $M_0(\varepsilon)$ and the matrices 
    	$e^{2 \pi\,i\,J_j(\varepsilon)}, j=L, R$ satisfy the following relation
    	$$\,
    	 M_0(\varepsilon)=e^{2 \pi\,i\,J_L(\varepsilon)}\,e^{2 \pi\,i\,J_R(\varepsilon)}=
    	 e^{2 \pi\,i\,J_R(\varepsilon)}\,e^{2 \pi\,i\,J_L(\varepsilon)}\,.
    	\,$$
    	\eco

    	\proof
    	
    	The statement follows immediately from \eqref{M-hat} and \eqref{M-0-B}.
    	
   	\qed

     Combining  \coref{J}  and \thref{an-mon} we have that
     
       \bpr{main1}
         Assume that the condition \eqref{res1} holds. Then the matrices $J_k(\varepsilon), k=L, R$
         of the unfolded equation and the  monodromy matrix $M_0$ of the DCHE are connected
         by the limit
          $$\,
           e^{2 \pi\,i\,J_L(\varepsilon)}\,e^{2 \pi\,i\,J_R(\varepsilon)}=
             e^{2 \pi\,i\,J_R(\varepsilon)}\,e^{2 \pi\,i\,J_L(\varepsilon)}
             \longrightarrow M_0
          \,$$
         when $\sqrt{\varepsilon} \rightarrow 0$.
       \epr
       
       Thanks to \prref{main1} we find the limit of $q_R+q_L$ when $\sqrt{\varepsilon} \to 0$.
       
       \bco{q}
        Assume that the condition \eqref{res1} holds. Then
          \be\label{limit}
           \lim_{\sqrt{\varepsilon} \to 0} q_R + q_L=\lambda\,,
          \ee
          where $\lambda$ is given by \eqref{a}.
       \eco
       
  In the Appendix we demonstrate by a direct computation that the limit \eqref{limit} is valid for lower values of the 
  parameter $\alpha$.

%%%%%%%%%%%%%%%%%%%%%%%%%%%%%%%%%%%%%%%%%%%%%%%%%%%%%%%%%%%%%%%%%%%%%%%%%%%%%%%
%stokes
%%%%%%%%%%%%%%%%%%%%%%%%%%%%%%%%%%%%%%%%%%%%%%%%%%%%%%%%%%%%%%%%%%%%%%%%%%         	

     \section{Unfolding of the Stokes matrix at the origin }
        
        Throughout this section we assume that $\beta \neq 0$ and therefore the origin is an irregular singularity for the DCHE.
           In \cite{St2} we have shown by a direct computation, that when $\alpha=2$  
        the Stokes matrix at the origin of the DCHE \eqref{eq1} can be obtained by
        a limit of this  part of the monodromy matrices around resonant singular points
         that governs the existence of the
        logarithmic term in the solution of the unfolded equation.
        In the  section 4.1 we  show that this result remains valid for every reducible DCHE \eqref{eq1}. 
        In fact the realization of the Stokes matrix as a limit of the part of the monodromy matrix of the unfolded equation
        is an effect of the recent theoretical result of Hurtubise, Klime\v{s}, Lambert and Rousseau \cite{HLR, Kl-2, CL-CR, CL-CR1}.
        Using the obtained connection between the analytic invariants of both equations we provide the
        Stokes matrix at the origin without studying in details the DCHE. Instead we deal with the unfolded
        equation and its monodromy matrix around a resonant singularity $x_L$.
        In the section 4.2 we build explicitly an actual fundamental matrix solution of
         the DCHE \eqref{eq1} at the origin with respect to which the Stokes matrix has the form obtained in the firs part.

        \subsection{The Stokes matrix at the origin as a limit of the monodromy matrix around $x_L$} 
       Following Lambert and Rousseau \cite{CL-CR, CL-CR1} we consider both equations in the ramified domain
       $\{x\in\CC\,:\, -\kappa < \arg (x) < \kappa\}$ where $0 < \kappa < \frac{\pi}{2}$. We cover this domain by two open sectors
         \ben
             \Omega_1=\Omega_1(\rho, \kappa) &=&
             \left\{ x=r\,e^{i \delta}\,|\,
            0 < r < \rho,\, -\kappa-\pi < \delta < \kappa \right\},\\[0.1ex]
              \Omega_2=\Omega_2(\rho, \kappa) &=&
              \left\{ x=r\,e^{i \delta}\,|\,
              0 < r < \rho,\, -\kappa < \delta < \kappa + \pi \right\}\,.
         \een
        The radius $\rho$ is so chosen that 
         $x_{LL}, x_{RR} \notin \Omega_1 \cup \Omega_2$ while $x_L, x_R\in \Omega_1 \cup \Omega_2$. 
         Denote by $\Omega_R$ and $\Omega_L$ the connected components of the
         intersection $\Omega_1\cap \Omega_2$, as $x_R\in \Omega_R,\,x_L\in \Omega_L$. 
         We have a proposition similar to \prref{con}. 
        
        \bpr{con-1}
        When $\sqrt{\varepsilon} \to 0$ the fundamental set of solutions $\{w_1(x, \varepsilon), w_2(x, \varepsilon)\}$
        of the unfolded equation fixed by \thref{t1-p} converges uniformly on compact sets of $\Omega_R \cup \Omega_L$ to the fundamental set of solutions
        $\{w_1(x, 0), w_2(x, 0)\}$ of the DCHE fixed by \thref{t1}.
        \epr

           \begin{figure}[t]
           	
           	\begin{minipage}[c][1\width]
           		{0.003\textwidth}
           		\centering
           		\includegraphics[scale=0.38]{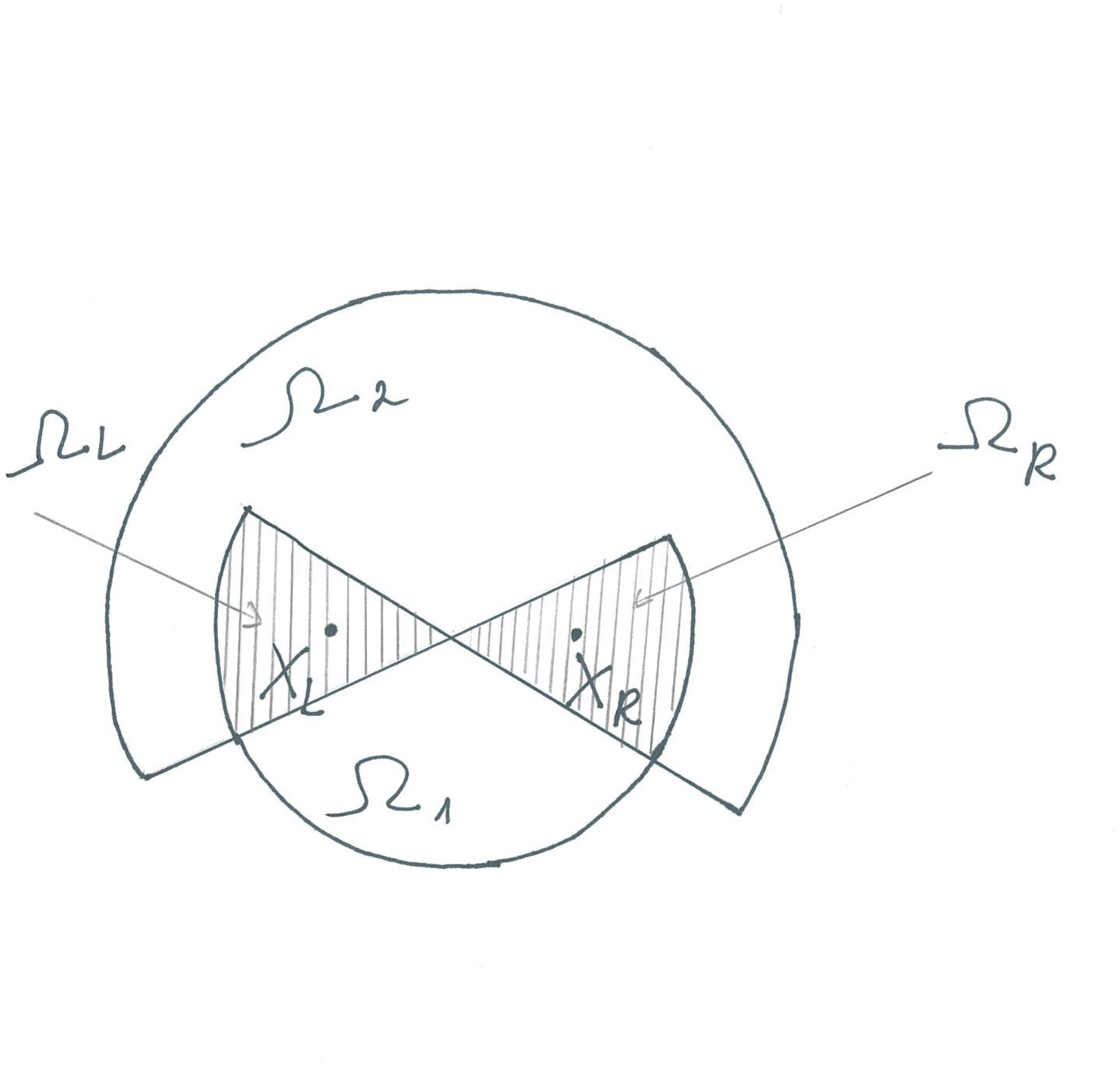}
           	\end{minipage}
           	\hfill
           	\begin{minipage}[c][1\width]
           		{0.5\textwidth}
           		\centering
           		\includegraphics[scale=0.38]{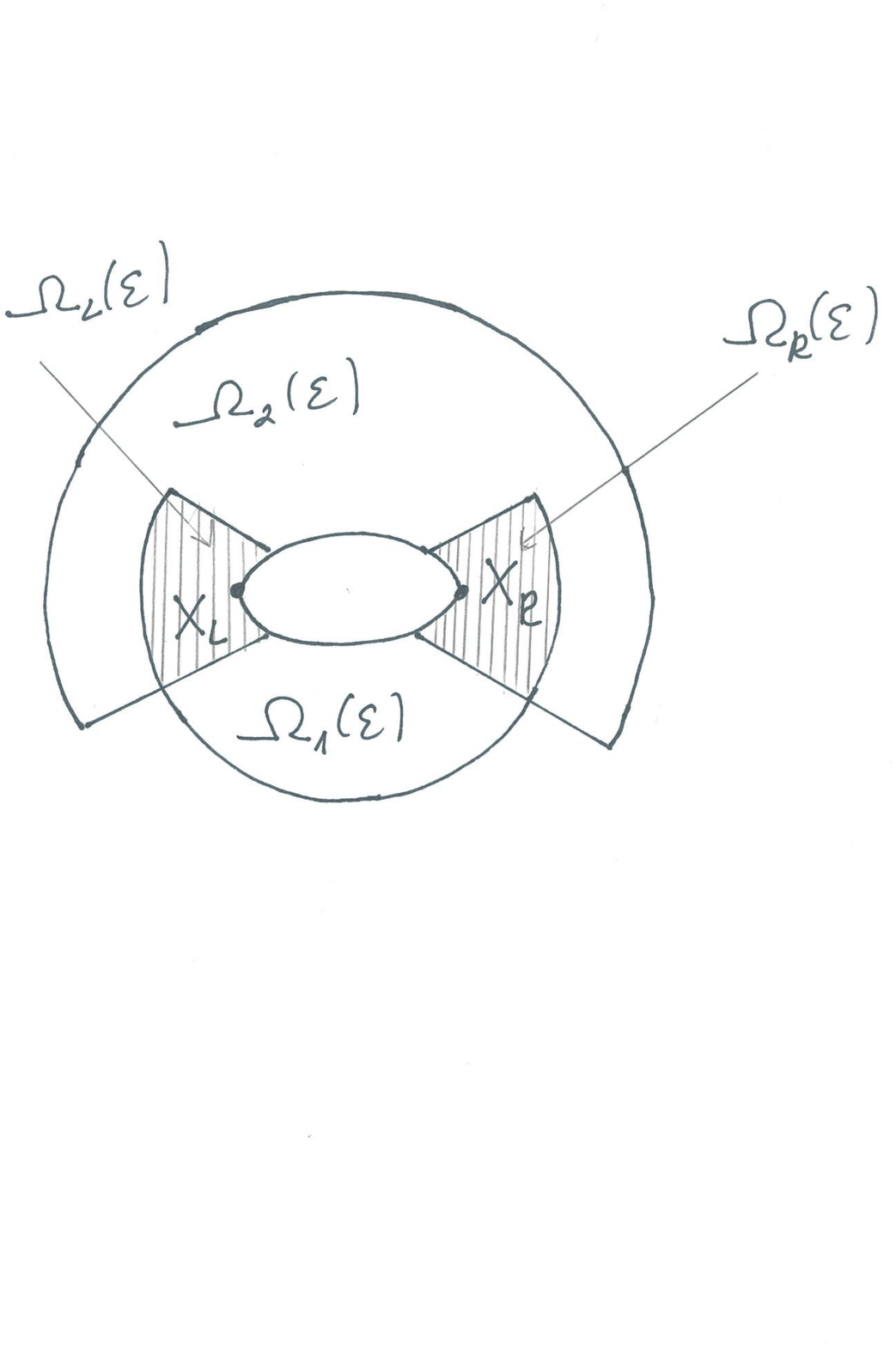}
           	\end{minipage}
           	\caption{The sectors $\Omega_j$ and $\Omega_j(\varepsilon),\,j=1, 2$.}
           	\label{fig:sectorimage}
           \end{figure}

      Consider the DCHE \eqref{eq1} over $\Omega_1 \cup \Omega_2$.          
     From the sectorial normalization theorem of Sibuya \cite{S} and the theorem of Hukuhara-Turrittin\cite{W} it follows that
     the actual fundamental matrix solution $\Phi_0(x, 0)$ of the DCHE can be represented as
      $$\,
        \Phi_j(x, 0)=\exp (G x)\,H_j(x)\,\left[x^{\Lambda}\,\exp\left(-\frac{B}{x}\right)\right]_j
      \,$$ 
      on the sectors $\Omega_j,\,j=1, 2$, respectively. Here
        \be\label{GLB}
        G=\diag (0, \gamma),\quad \Lambda=\diag(0, -\alpha), \quad B=\diag(0, \beta)
                \ee
       and $[x^{\Lambda}\,\exp(-B/x)]_j$ is the branch of the matrix $x^{\Lambda}\,\exp(-B/x)$ on $\Omega_j, j=L, R$,
       respectively.
       The matrices $H_j(x)$ are  holomorphic matrix functions on $\Omega_j$, respectively, as both of them are asymptotic
       in the Gevrey 1 sense to the same formal matrix $\hat{H}(x)$ on $\Omega_j, j=L, R$. On the sector $\Omega_R$
       the fundamental matrix solutions $\Phi_1(x, 0)$ and $\Phi_2(x, 0)$ coincide.
       On the sector $\Omega_L$ the jump of the solution $\Phi_2(x, 0)$ to the solution $\Phi_1(x, 0)$ is measured geometrically
       by the  Stokes matrix $St_{\pi}$, corresponding to the singular direction $\theta=\pi$
          \be\label{s-m}
           \Phi_2(x, 0)=\Phi_1(x, 0)\,St_{\pi}\,,
          \ee   
       where
        $$\,
               St_{\pi}=\left(\begin{array}{cc}
                 1   & \mu\\
                 0   & 1
                            \end{array}\right)\,.
             \,$$ 

    Consider now the unfolded equation. The next theorem describes the behavior of the fundamental matrix solution $\Phi_0(x, \varepsilon)$
    near the singular points when $x_L$ is a resonant singularity.

   \bth{f-p}
          Assume that the condition \eqref{res2} holds. Then 
          the fundamental matrix solution $\Phi_0(x, \varepsilon)$ of the unfolded equation depends analytically on $\sqrt{\varepsilon}$ and it is represented
             in a neighborhood of the resonant singularity $x_L$ which does not contain the point $x_R$ as
                \be\label{r-F}
                 \Phi_0(x, \varepsilon)=\left(I_L(\varepsilon) + \mathcal{O}(x-x_L)\right)\,
                 (x-x_L)^{\frac{1}{2} \Lambda + \frac{1}{2 x_L}\,B}\,(x-x_L)^{T_L(\varepsilon)}
                \ee
                and in a neighborhood of the non-resonant singularity $x_R$ which does not contain the point $x_L$ as
                    $$\,
                    \Phi_0(x, \varepsilon)=\left(I_R(\varepsilon) + \mathcal{O}(x-x_R)\right)\,
                    (x-x_R)^{\frac{1}{2} \Lambda + \frac{1}{2 x_R}\,B}\,.
                  \,$$
                The matrices $I_j(\varepsilon) + \mathcal{O}(x-x_j),\,j=L, R$ are  holomorphic matrix functions there.
                The matrices $\Lambda$ and $B$ are given by \eqref{GLB}. The matrix $T_L(\varepsilon)$ is defined as
                 $$\,
                  T_L(\varepsilon)=\left(\begin{array}{cc}
                    0   & d_L\\
                    0   & 0
                            \end{array}\right)\,,
                 \,$$ 
                 where
                 $$\,
                  d_L=\res (w_2'(x, \varepsilon),\,x=x_L)\,.
                 \,$$
                 \ethe
       
         \proof
         The proof is similar to the proof of Proposition 4.7 in \cite{St1}.
         \qed

        With respect to the matrix solution $\Phi_0(x, \varepsilon)$ from \thref{f-p} the local monodromy matrices
        $M_j(\varepsilon)$
       around the singular point $x_j, j=L, R$ have the form
         $$\,
          M_L(\varepsilon)=e^{2\pi\,i\,(\frac{1}{2}\,\Lambda + \frac{1}{2\,x_L}\,B)}\,e^{2 \pi\,i\,T_L(\varepsilon)}=
          e^{2 \,\pi\,T_L(\varepsilon)},\qquad
           M_R(\varepsilon)=e^{2\pi\,i\,(\frac{1}{2}\,\Lambda + \frac{1}{2\,x_R}\,B)}\,.
         \,$$
         
        Let $\Omega_1(\varepsilon)$ and $\Omega_2(\varepsilon)$ be the sectors obtained from the sectors $\Omega_1$ and $\Omega_2$
        by making a cut between the points $x_L$ and $x_R$ through the real axis (see  Figure \ref{fig:sectorimage}). 
        The origin belongs to this cut. When 
        $\sqrt{\varepsilon} \to 0$ the sectors $\Omega_j(\varepsilon)$ tend to the sectors $\Omega_j, j=L, R$, respectively.
         Consider the unfolded equation over $\Omega_1(\varepsilon) \cup \Omega_2(\varepsilon)$. 
        The fundamental matrix solution $\Phi_0(x, \varepsilon)$ writes  also as
        \be\label{r-F-p}
         \Phi_0(x, \varepsilon)=G(x, \varepsilon)\,H(x, \varepsilon)\,
         (x-x_L)^{\frac{1}{2}\,\Lambda + \frac{1}{2\, x_L}\,B}\,
          (x-x_R)^{\frac{1}{2}\,\Lambda + \frac{1}{2\, x_R}\,B}\,,
        \ee  
        where
         $$\,
          G(x, \varepsilon)=(x-x_{LL})^{-\frac{x_{LL}}{2}\,G}\,(x_{RR} -x)^{-\frac{x_{RR}}{2}\,G}\,
         \,$$
         and
         $$\,
          H(x, \varepsilon)=\left(\begin{array}{cc}
             1     & (x-x_L)^{-\frac{\beta}{2 \sqrt{\varepsilon}}+ \frac{\alpha}{2}}\,
                   (x-x_R)^{\frac{\beta}{2 \sqrt{\varepsilon}}+ \frac{\alpha}{2}}\,w_2(x, \varepsilon)\\[0.1ex]
             0     & 1        
                       \end{array}
                       \right)\,.
         \,$$
       The next proposition that follows immediately from \prref{con-1} is a key for the extension of the results in \cite{St2}
       
       \bco{im}
         The fundamental matrix solution $\Phi_0(x, \varepsilon)$ from \eqref{r-F-p} converges uniformly on compact sets of
         $\Omega_R \cup \Omega_L$ to the actual fundamental matrix solution $\Phi_0(x, 0)$ at the origin of the DCHE when
         $\sqrt{\varepsilon} \to 0$.
       \eco

         \begin{figure}[h]
         	
         	\centering
         	
         	\includegraphics[scale=0.3]{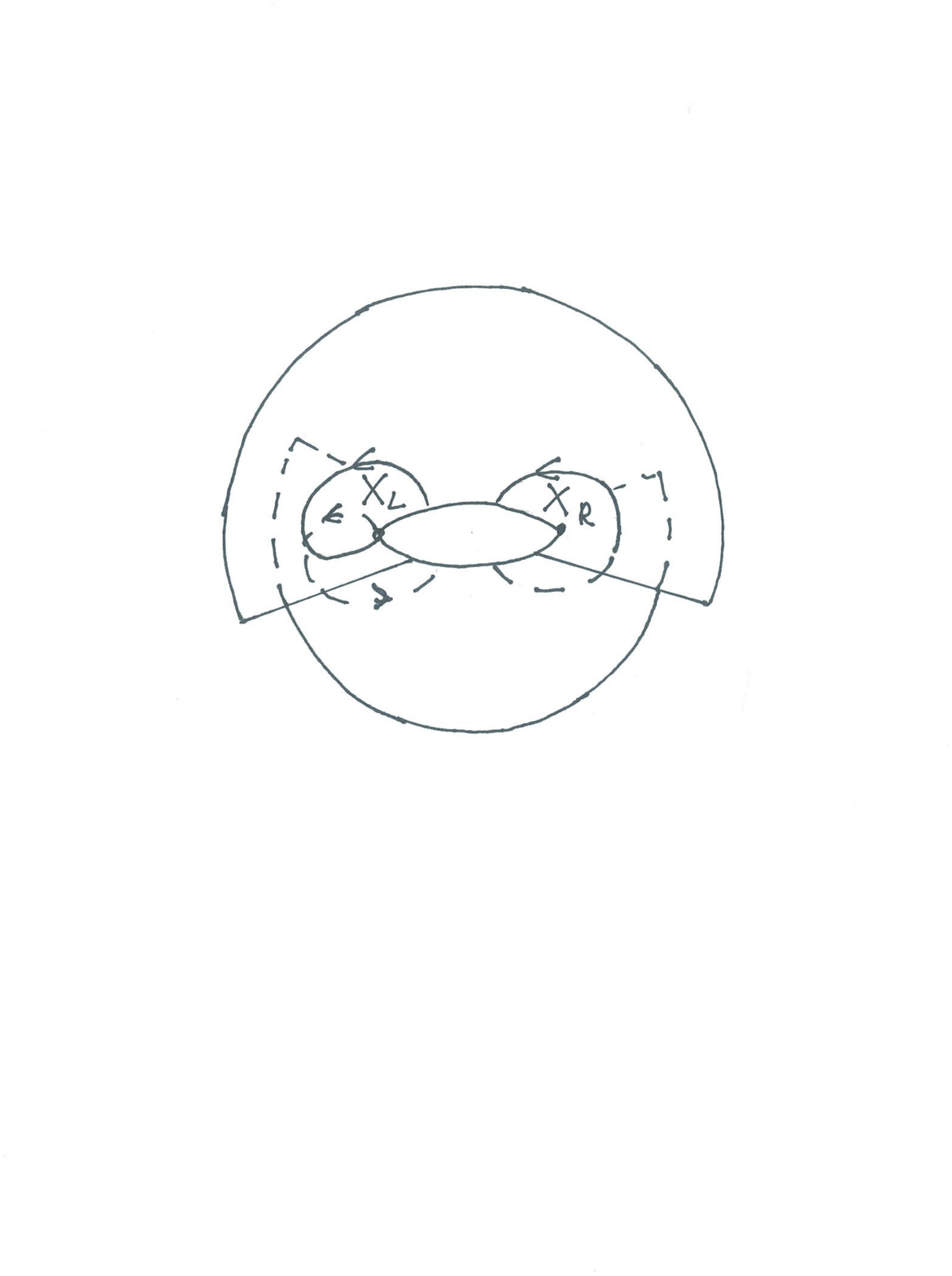}
         	\caption{The monodromy operator  $M_L(\varepsilon)$.}
         	\label{fig:stokes}
         \end{figure}

       As in \cite{St2} \coref{im} allows us to identify the so called unfolded Stokes matrix $St_L(\varepsilon)$ with the
       matrix $e^{2 \pi\,i\,T_L}$ when the point $x_L$ is a resonant singularity.
       
       \bpr{m-s}$($ Proposition 6.1 in \cite{St2}$)$
       Let $M_L(\varepsilon)$ and $St_L(\varepsilon)$ be the monodromy matrices and the unfolded Stokes matrix of the
       unfolded equation. Then when the condition  \eqref{res2} holds they satisfy the relations
        $$\,
          M_L(\varepsilon)=St_L(\varepsilon)\,e^{\pi\,i (\Lambda + \frac{1}{x_R}\,B)}
        \,$$
        on the sector $\Omega_1(\varepsilon)$, and
          $$\,
          M_L(\varepsilon)=e^{\pi\,i (\Lambda + \frac{1}{x_R}\,B)}\,St_L(\varepsilon)
          \,$$
          on the sector $\Omega_2(\varepsilon)$.
          \epr

       \bco{c1}$($Corollary 6.2 in \cite{St2}$)$
        Assume that the condition \eqref{res2} holds. Then
           $$\,
             St_L(\varepsilon)=e^{2\,\pi\,i\,T_L(\varepsilon)}\,.
           \,$$
       \eco
       Let us compute the limit of the matrix $T_L(\varepsilon)$ when $\sqrt{\varepsilon} \to 0$.
              The next proposition and lemma  gives us the number $d_L$and its limit when $\sqrt{\varepsilon} \to 0$.
         
       \bpr{d}
        Assume that the condition \eqref{res2} holds. Then 
         \ben
          d_L  &=&
          \frac{1}{\left(\frac{\beta}{2 \sqrt{\varepsilon}}+\frac{\alpha}{2}-1\right)!}
                   \sum_{k=0}^{\frac{\beta}{2 \sqrt{\varepsilon}} + \frac{\alpha}{2}-1}
            \left(\begin{array}{c}
               \frac{\beta}{2 \sqrt{\varepsilon}} + \frac{\alpha}{2} -1\\[0.1ex]
                 k
                 \end{array}\right)\,
                \frac{\Gamma \left(\frac{\beta}{2 \sqrt{\varepsilon}}-\frac{\alpha}{2}+1\right)(-2 \sqrt{\varepsilon})^{-\alpha+1+k}}
                {\Gamma (2+k-\alpha)}\,A_L\,,
         \een
         where $A_L$ is given by \eqref{al}.
          \epr

       \ble{d}
             Assume that the condition \eqref{res2} holds. Then 
             $$\,
              \lim_{\sqrt{\varepsilon} \to 0} d_L=
              (-\beta)^{1-\alpha}\,\sum_{k=0}^{\infty} 
              (-1)^k\,\frac{\beta^k\,\gamma^k}{k!\,\Gamma(2+k-\alpha)}\,.
             \,$$
       \ele
       \proof
       Applying the limit \eqref{l} for $z=\frac{\gamma}{2 \sqrt{\varepsilon}}$ we find that $A \to \gamma^k$
       when $\sqrt{\varepsilon} \to 0$. Again applying the limit \eqref{l} for $z=\frac{\beta}{2 \sqrt{\varepsilon}}$ and using the limit
       \eqref{ll} we obtain the limit of $d_L$.
       \qed
       
	  In \cite{CL-CR1} Lambert and Rousseau prove that the unfolded Stokes matrix $St_L(\varepsilon)$ depends 
	  analytically on $\sqrt{\varepsilon}$ and tends to the Stokes matrix $St_{\pi}$ when $\sqrt{\varepsilon} \to 0$.
	  Then
	  
	  \bth{limit}
	   Assume that the condition \eqref{res2} holds. Then the Stokes matrix $St_{\pi}$ at the origin of the DCHE
	   and the matrix $e^{2 \pi\,i\,T_L(\varepsilon)}$ of the unfolded equation are connected as
	    $$\,
	      e^{2 \pi\,i\,T_L(\varepsilon)} \longrightarrow St_{\pi}
	    \,$$
	    when $\sqrt{\varepsilon} \to 0$.
	  \ethe
	  
	  In a consequence of \thref{limit} we have
	   
	\bco{stokes}
	 Assume that $\beta > 0 $. Then there exists an actual fundamental matrix solution at the origin of the DCHE
	 with respect to which the corresponding Stokes matrix $St_{\pi}$ is given by
	  $$\,
	    St_{\pi}=\left(\begin{array}{cc}
	       1   & \mu\\
	       0   & 1
	       \end{array}\right)\,,
	  \,$$
	  where
	  \be\label{sm}
	   \mu= 2\,\pi\,i\,  (-\beta)^{1-\alpha}\,\sum_{k=0}^{\infty} 
	   (-1)^k\,\frac{\beta^k\,\gamma^k}{k!\,\Gamma(2+k-\alpha)}\,.
	  \ee
	\eco

	\subsection{Actual fundamental matrix solution at the origin of the DCHE }
	
	In this paragraph we present explicitly the actual fundamental matrix solution $\Phi_0(x, 0)$ at the origin  of the DCHE
	with respect to which the Stokes matrix $St_{\pi}$ has the form fixed by \thref{stokes}. 
	When $\alpha\in\ZZ$ we apply the Borel-Laplace summation in order to build this actual solution (see \cite{LR, R2} for details).
	When $\alpha\notin\ZZ$ we directly express the solution $w_2(x, 0)$ in terms of Laplace integrals without using the summability 
	theory. The second approach is more general and also can be used when $\alpha\in\ZZ$. But the first approach allows us to
	distinguish special solutions of DCHE that are holomorphic in whole $\CC^*$ even when $\alpha\notin\ZZ$.
	
	We start by building a formal fundamental matrix solution at the origin.
	 \bth{0}
	 Assume that $\beta >  0$. Then the DCHE \eqref{eq1}  possesses an unique formal fundamental matrix solution $\hat{\Phi}_0(x, 0)$
	 at the origin in the form
	 $$\,
	   \hat{\Phi}_0(x, 0)=\exp(G x)\,\hat{H}(x)\,x^{\Lambda}\,\exp \left(-\frac{B}{x}\right)\,,
	 \,$$
	 where the matrices $G, \Lambda$ and $B$ are given by \eqref{GLB}.
	 The matrix $\hat{H}(x)$ is defined as follows:
	 \begin{enumerate}
	 	
	 	\item\,
	 	If $\alpha\notin \ZZ$ then
	 	  $$\,
	 	  \hat{H}(x)=\left(\begin{array}{cc}
	 	    1   & \frac{x\,\hat{\varphi}(x)}{\beta}\\[0.1ex]
	 	    0   & 1
	 	    \end{array}\right)
	 	  \,$$
	 	  where
	 	  \be\label{s1}
	 	   \hat{\varphi}(x)=\sum_{n=0}^{\infty}
	 	    \frac{(-1)^n\,S_n\,\Gamma(2-\alpha+n)}{\beta^n}\,x^{n+1}\,.
	 	  \ee
	 	  Here $S_n$ is the $n$-th partial sum of the absolutely convergent number series
	 	    $$\,
	 	      S=\sum_{k=0}^{\infty}
	 	       \frac{(-1)^k\,\beta^k\,\gamma^k}{k!\,\Gamma(2-\alpha+k)}\,.
	 	    \,$$
	 	    In particular,
	 	    \begin{enumerate}
	 	    	\item\,
	 	    	If $S=0$ then the power series $\hat{\varphi}(x)$ is convergent.
	 	    	
	 	    	\item\,
	 	    	If $S \neq 0$ then the power series $\hat{\varphi}(x)$ is divergent.
	 	    \end{enumerate}

	 	\item\,
	 	If $\alpha\in\NN$ then  	
	 	 	  $$\,
	 	 	  \hat{H}(x)=\left(\begin{array}{cc}
	 	 	  1   & \frac{x^{\alpha}\,\gamma^{\alpha-1}\,\hat{\psi}(x)}{\beta} + x^{\alpha}\,P\left(\frac{1}{x}\right)\\[0.1ex]
	 	 	  0   & 1
	 	 	  \end{array}\right)
	 	 	  \,$$
	 	 	  where
	 	 \be\label{s2}
	 	\hat{\psi}(x)=
	 \sum_{n=0}^{\infty}
	 	\frac{(-1)^n\,W_n\,n!}{\beta^n}\,x^{n+1}
	 	\ee
	 	and $P\left(\frac{1}{x}\right)$ is a polynomial in $\frac{1}{x}$ of degree $\alpha-2$ for $\alpha \geq 2$ and
	 	$P\equiv 0$ for $\alpha=1$.
	 	Here $W_n$ is the $n$-th partial sum of the absolutely convergent number series
	 	$$\,
	 	W=\sum_{k=0}^{\infty}
	 	\frac{(-1)^k\,\beta^k\,\gamma^k}{k!\,\Gamma(\alpha+k)}\,.
	 	\,$$
	 	In particular,
	 	\begin{enumerate}
	 		\item\,
	 		If $W=0$ then the power series $\hat{\psi}(x)$ is convergent.
	 		
	 		\item\,
	 		If $W\neq 0$ then the power series $\hat{\psi}(x)$ is divergent.
	 	\end{enumerate}
	 	
	 	\item\,
	 	If $\alpha\in\ZZ_{\leq 0}$ then 
	 	  	  $$\,
	 	  	  \hat{H}(x)=\left(\begin{array}{cc}
	 	  	  1   & \frac{x\,\hat{\phi}(x)}{\beta}\\[0.1ex]
	 	  	  0   & 1
	 	  	  \end{array}\right)
	 	  	  \,$$
	 	  	  where
	 	\be\label{s3}
	 	\hat{\phi}(x)=\sum_{n=0}^{\infty}
	 	\frac{(-1)^n\,Q_n\,\Gamma(2-\alpha+n)}{\beta^{1-\alpha+n}}\,x^{n+1}\,.
	 	\ee
	 	Here $Q_n$ is the $n$-th partial sum of the absolutely convergent number series
	 	$$\,
	 	Q=\sum_{k=0}^{\infty}
	 	\frac{(-1)^k\,\gamma^k\,\beta^{1-\alpha+k}}
	 	{k!\,\Gamma(2-\alpha+k)}\,.
	 	\,$$    
	 	In particular,
	 	\begin{enumerate}
	 		\item\,
	 		If $Q=0$ then the power series $\hat{\phi}(x)$ is convergent.
	 		
	 		\item\,
	 		If $Q\neq 0$ then the power series $\hat{\phi}(x)$ is divergent.
	 	\end{enumerate}
	 \end{enumerate}
	 \ethe
	 
    \proof
    
    The proof is similar to the proof of Proposition 4.2 in \cite{St2}. We only note that when $\alpha\notin\ZZ$ we reduce the solution
    $w_2(x, 0)$ from \eqref{fss} to the integral $\int_0^x \frac{e^{-\frac{\beta}{z}}}{z^{\alpha}}\, d z$ while when $\alpha\in\ZZ$
    to the integral $\int_0^x \frac{e^{-\frac{\beta}{z}}}{z}\,d z$.
    \qed

	The application of summability theory
	 to the differential equations ensures that the divergent power series $\hat{\varphi}(x), \hat{\psi}(x)$ and $\hat{\phi}(x)$ are 1-summable
	 in any direction $\theta$  except for the singular direction $\theta=\pi$.
	 
	 \ble{sum}
	 1.\,Assume that $W \neq 0$. Then for any direction $\theta \neq \pi$ the function
	 $$\,
	 \psi_{\theta}(x)=\beta\,\int_0^{+\infty\,e^{i \theta}}
	 \frac{v(\xi)\,e^{-\frac{\xi}{x}}}{\xi+ \beta}\,d \xi
	 \,$$
	 defines the 1-sum of the power series $\hat{\psi}(x)$ from \eqref{s2} in such a direction. Here 
	 $$\,
	 v(\xi)=\sum_{n=0}^{\infty} \frac{\gamma^n\,\xi^n}{n!\,\Gamma(\alpha+n)}\,.
	 \,$$

	 2.\,  Assume that $Q \neq 0$. Then for every direction $\theta \neq \pi$ the  function
	 $$\,
	 \phi_{\theta}(x)=
	 \frac{\beta}{x^{1-\alpha}}\int_0^{+\infty e^{i \theta}}
	 \frac{\xi^{1-\alpha}\,q(\xi)\,e^{-\frac{\xi}{x}}}{\xi+\beta}\,d \xi
	 \,$$
	 defines the 1-sum of the  power series $\hat{\phi}(x)$ from \eqref{s3} is such a direction.
	 Here
	 $$\,
	 q(\xi)=\sum_{k=0}^{\infty} \frac{\gamma^k\,\xi^k}{k!\,\Gamma(2-\alpha+k)}\,.
	\,$$
	 
	 The functions  $\psi_{\theta}(x)$ and $\phi_{\theta}(x)$ are holomorphic functions in the  open disc
	    \be\label{D}
	    \mathcal{D}_{\theta}(|\gamma|)=\left\{
	    x\in\CC^*\,|\, Re \,\left(\frac{e^{i \theta}}{x}\right) > |\gamma|\right\}\,.
	    \ee
	 \ele

	 \proof
	 We will prove the second statement of the Lemma. The first is proved in a similar way.
	 
	 Let $\alpha\in\ZZ_{\leq 0}$. The formal Borel transform of order 1 of the series $\hat{\phi}(x)$ from
	 \eqref{s3} yields the convergent power series near the origin $\xi=0$ of the Borel $\xi$-plane
	 $$\,
	 \phi(\xi)=(\hat{\B}_1 \hat{\phi})(\xi)=
	 \sum_{n=0}^{\infty} (-1)^n\,Q_n\,\frac{(1-\alpha+n)!}{\beta^{1-\alpha+n}}\,\frac{\xi^n}{n!}\,.
	 \,$$
	 The series $\phi(\xi)$ can be regarded as the $1-\alpha$-th derivative of the series
	 $$\,
	 w(\xi)=\sum_{n=0}^{\infty} (-1)^n\,Q_n\,\frac{\xi^{1-\alpha+n}}{\beta^{1-\alpha+n}}\,.
	 \,$$
	 Then from Lemma 6.12 in \cite{St2} it follows that the Laplace transform of order 1 of the function
	 $\phi(\xi)$ is expressed by the Laplace transform of the function $w(\xi)$
	 $$\,
	 (\L_1 \phi)(x)=\frac{1}{x^{1-\alpha}}\,(\L_1 w)(x)
	 \,$$
	 since $\frac{d^k\,w}{d\, \xi^k}=0$ for $0 \leq k \leq -\alpha$. It turns out that the power series
	 $w(\xi)$ is the Maclaurin series of a well known function. More precisely,
	 $$\,
	 w(\xi)=(-\beta)^{1-\alpha}\,
	 \left(\frac{1}{1+\frac{\xi}{\beta}} - \sum_{k=0}^{-\alpha} (-1)^k\,\frac{\xi^k}{\beta^k}\right)\,
	 q(\xi)\,,
	 \,$$ 
	 where $q(\xi)$ is the analytic function in $\CC$ defined by the power series
	 $$\,
	 q(\xi)=\sum_{k=0}^{\infty} \frac{\gamma^k\,\xi^k}{k!\,\Gamma(2-\alpha+k)}\,.
	 \,$$
	 Since
	 $$\,
	 \sum_{k=0}^{-\alpha} (-1)^k\frac{\xi^k}{\beta^k}=
	 \frac{1-(-1)^{1-\alpha}\frac{\xi^{1-\alpha}}{\beta^{1-\alpha}}}{1+\frac{\xi}{\beta}}
	 \,$$
	 then
	 \ben
	 (\L_1 \phi)(x)  &=&
	 \frac{(-\beta)^{1-\alpha}}{x^{1-\alpha}}\,
	 \left(\L_1 \frac{q(\xi)}{1+\frac{\xi}{\beta}}\right)(x) -
	 \frac{(-\beta)^{1-\alpha}}{x^{1-\alpha}}\,
	 \left(\L_1 \frac{q(\xi)}{1+\frac{\xi}{\beta}}\right)(x) \\[0.3ex]
	 &+&
	 \frac{1}{x^{1-\alpha}} \left(\L_1 \frac{\xi^{1-\alpha}\,q(\xi)}{1+\frac{\xi}{\beta}}\right)(x)=
	 \frac{\beta}{x^{1-\alpha}}\int_0^{+\infty e^{i \theta}}
	 \frac{\xi^{1-\alpha}\,q(\xi)\,e^{-\frac{\xi}{x}}}{\xi+\beta}\,d \xi\,.    
	 \een
	 Thus the function
	 $$\,
	 \phi_{\theta}(x)=
	 \frac{\beta}{x^{1-\alpha}}\int_0^{+\infty e^{i \theta}}
	 \frac{\xi^{1-\alpha}\,q(\xi)\,e^{-\frac{\xi}{x}}}{\xi+\beta}\,d \xi
	 \,$$
	 gives the 1-sum of the power series $\hat{\phi}(x)$ in any direction $\theta \neq \pi$.
    Since 
      $$\,
      \left|\frac{q(\xi)}{\xi+\beta}\right| \leq A \sum_{k=0}^{\infty} \frac{|\gamma|^k\,|\xi|^k}{k!}=
     A\,e^{|\gamma|\,|\xi|}
      \,$$ 
      for an appropriate constant $A > 0$ the integral $\phi_{\theta}(x)$ exists when $\alpha\in\ZZ_{\leq 0}$
    and defines a holomorphic function in the open disc $Re\,\left(\frac{e^{i \theta}}{x}\right) > |\gamma|$.
    
    This ends the proof.
	 \qed
	
	\bre{s1}
	 Unfortunately till now we can not derive the 1-sum of the power series $\hat{\varphi}(x)$ from
	 \eqref{s1} in an explicit way. The formal Borel transform of this series
	  $$\,
	   (\hat{\B}_1 \hat{\varphi})(\xi)=
	   \sum_{n=0}^{\infty} \frac{(-1)^n\,S_n\,\Gamma(2-\alpha+n)}{\beta^n}\,
	     \frac{\xi^n}{n!}
	  \,$$
	  is a convergent power series for $|\xi| < \beta$. But we can not specify explicitly the function
	  whose Maclaurin series is $(\hat{\B}_1 \hat{\varphi})(\xi)$. For this reason we use a slightly different
	  approach to build an actual solution of the DCHE when $\alpha\notin\ZZ$. 
	\ere

	\bre{ex}
			Let $I=(-\pi, \pi) \subset \RR$.
	 When we move the direction $\theta\in I$ the holomorphic functions 
	$\psi_{\theta}(x)$ (resp. $\phi_{\theta}(x)$) glue together analytically and define a holomorphic function $\tilde{\psi}(x)$ 
	(resp. $\tilde{\phi}(x)$) on a sector $\widetilde{\mathcal{D}}$
	with opening $> \pi$
	 \be\label{D}
	  \widetilde{\mathcal{D}}=\bigcup_{\theta\in I}\,\widetilde{\mathcal{D}}_{\theta}(|\gamma|)\,,
	 \ee
	 where $\widetilde{\mathcal{D}}_{\theta}(|\gamma|)$ is the lifting of $\mathcal{D}_{\theta}(|\gamma|)$
	 on the Riemann surface of the natural logarithm. On $\widetilde{\mathcal{D}}$ the function $\tilde{\psi}(x)$ (resp. $\tilde{\phi}(x)$) 
	 is asymptotic to
	 the power series $\hat{\psi}(x)$ (resp. $\hat{\phi}(x)$) in Gevrey 1 sense and defines the 1-sum of this series there. The restriction of
	 $\tilde{\psi}(x)$ (resp. $\tilde{\phi}(x)$) on $\CC^*$ is a multivalued function. In every direction $\theta\neq \pi$ the function $\tilde{\psi}(x)$ (resp. $\tilde{\phi}(x)$)
	 has only one value than coincides with the function $\psi_{\theta}(x)$ from \leref{sum}(1) (resp. $\phi_{\theta}(x)$ from \leref{sum}(2)). Near the singular direction
	 $\theta=\pi$ the function $\tilde{\psi}(x)$ (resp. $\tilde{\phi}(x)$) has two different values: $\psi^{+}_{\pi}(x)=\psi_{\pi+\epsilon}(x)$  (resp. $\phi^{+}_{\pi}(x)=\phi_{\pi+\epsilon}(x)$) and
	 $\psi^{-}_{\pi}(x)=\psi_{\pi-\epsilon}(x)$ (resp. $\phi^{-}_{\pi}(x)=\phi_{\pi-\epsilon}(x)$) for a small number $\epsilon > 0$.
	\ere
	
	Now we can present an actual fundamental matrix solution at the origin.
	
	 \bth{0-1}
	 Assume that  $\beta >  0$. 
	 Then
	  \begin{enumerate}
	 
	  \item\,
	  Assume that one of the following conditions holds: $(\alpha\in\NN,\, W=0),\,(\alpha\in\ZZ_{\leq 0},\,Q=0)$
	  or $(\alpha\notin \ZZ, S=0)$. Then the DCHE \eqref{eq1} possesses an unique actual fundamental
	  matrix solution $\Phi_0(x, 0)$ at the origin in the form
	   \be\label{sol}
	    \Phi_0(x, 0)=\exp(G x)\,H(x)\,x^{\Lambda}\,\exp \left(-\frac{B}{x}\right)\,,
	   \ee
	   where the matrices $G, \Lambda$ and $B$ are given by \eqref{GLB} and $H(x)$ is a holomorphic matrix function in whole $\CC$.
	   More precisely, when $\alpha\notin\ZZ$ and $S=0$ the matrix $H(x)$ coincides with the matrix $\hat{H}(x)$ 
	   from \thref{0}(1). When
	   $\alpha\in\NN$ and $W=0$ the matrix $H(x)$ coincides with the matrix $\hat{H}(x)$ from \thref{0}(2).
	   When $\alpha\in\ZZ_{\leq 0}$ and $Q=0$ the matrix $H(x)$ coincides with the matrix $\hat{H}(x)$ from \thref{0}(3).
	 
	  \item\,
	  Assume that $Q\neq 0,\,W\neq 0, \,S\neq 0$. 
	  Then the DCHE \eqref{eq1} possesses an unique actual matrix solution $\widetilde{\Phi}_0(x, 0)$ at the origin in the form
	  \eqref{sol} which is a holomorphic matrix function on the sector $\widetilde{\mathcal{D}}$ from \eqref{D} whose opening is 
	  $>\pi$. The restriction of $\widetilde{\Phi}_0(x, 0)$ on $\CC^*$ is a multivalued function.
	  For any direction $\theta\neq \pi$   this solution has only one value  
	$\Phi^{\theta}_0(x, 0)$ in the form
	 \be\label{afm}
	 \Phi^{\theta}_0(x, 0)=\exp(G x)\,H_{\theta}(x)\,x^{\Lambda}\,\exp \left(-\frac{B}{x}\right)\,,
	 \ee
	 where the matrices $G, \Lambda$ and $B$ are given by \eqref{GLB}.
	 The matrix   $H_{\theta}(x)$ is defined as follows:
	 \begin{enumerate}
	 	
	 	\item\,
	 	If $\alpha\notin \ZZ$ then
	 	 $$\,
	 	 H_{\theta}(x)=\left(\begin{array}{cc}
	 	 1   & \frac{x\,\varphi_{\theta}(x)}{\beta}\\[0.15ex]
	 	 0   & 1
	 	 \end{array}\right)\,,
	 	 \,$$
	 	where
	 	\be\label{S}
	 	 \varphi_{\theta}(x)=\sum_{k=0}^{\infty}
	 	  \frac{\gamma^k\,x^k}{k!}\,
	 	  \left(\int_0^{+\infty e^{i \theta}}
	 	  \frac{e^{-\frac{\xi}{x}}}{(1+\frac{\xi}{\beta})^{2-\alpha+k}}\,d \xi\right)\,.
	 	\ee

	 \item\,
	 	If $\alpha\in\NN$ then  	
	 	$$\,
	 	H_{\theta}(x)=\left(\begin{array}{cc}
	 	1   & \frac{x^{\alpha}\,\gamma^{\alpha-1}\,\psi_{\theta}(x)}{\beta} + x^{\alpha}\,P\left(\frac{1}{x}\right)\\[0.1ex]
	 	0   & 1
	 	\end{array}\right)
	 	\,$$
	 	where $\psi_{\theta}(x)$ is defined by \leref{sum}(1) and extended by \reref{ex}.

	 	\item\,
	 	If $\alpha\in\ZZ_{\leq 0}$ then 
	 	$$\,
	    H_{\theta}(x)=\left(\begin{array}{cc}
	 	1   & \frac{x\,\phi_{\theta}(x)}{\beta}\\[0.1ex]
	 	0   & 1
	 	\end{array}\right)
	 	\,$$
	 	where $\phi_{\theta}(x)$ is defined by \leref{sum}(2) and extended by \reref{ex}.
	 	
	 	For the singular direction $\theta=\pi$ the DCHE \eqref{eq1} possesses two actual fundamental matrix solution at the
	 	origin
	 	  $$\,
	 	     \Phi_0^{\pi\pm}(x, 0)=\Phi_0^{\pi \pm \epsilon}(x, 0)\,,
	 	  \,$$
	 	  where $\epsilon > 0$ is a small number and the matrices $\Phi_0^{\pi \pm \epsilon}$ are given by \eqref{afm}.
	\end{enumerate}
	\end{enumerate}
	 \ethe

    \proof
    The proof of items $(1), (2.b)$ and $(2.c)$ follows directly from the theorem of Hukuhara-Turrittin-Martinet-Ramis \cite{JM-JR, R2}.
    We give the proof of item $(2.a)$.
    
    Let $\alpha\notin \ZZ$. The solution $w_2(x, 0)$ from \eqref{fss} becomes
     $$\,
      w_2(x, 0)=\sum_{k=0}^{\infty} \frac{\gamma^k}{k!}
      \left(\int_0^x \frac{e^{-\frac{\beta}{z}}}{z^{\alpha-k}}\,d z\right)=
        e^{-\frac{\beta}{x}}
        \sum_{k=0}^{\infty} \frac{\gamma^k}{k!}
        \left(\int_0^x \frac{e^{-\frac{\beta}{z}+ \frac{\beta}{x}}}{z^{\alpha-k}}\,d z\right)
        \,.     
      \,$$
    By setting
     $$\,
      - \frac{\beta}{z} + \frac{\beta}{x}=-\frac{\xi}{x}
     \,$$
    we transform  the solution $w_2(x, 0)$ into
     $$\,
      w_2(x, 0)=
      \frac{e^{-\frac{\beta}{x}}\,x^{1-\alpha}}{\beta}
      \sum_{k=0}^{\infty} \frac{\gamma^k\,x^k}{k!}
      \left(\int_0^{+\infty} 
      \frac{e^{-\frac{\xi}{x}}}{(1+\frac{\xi}{\beta})^{2-\alpha+k}}\,d \xi\right)
      \,.     
      \,$$
 
  Now we will show that this infinite sum defines a holomorphic function on the open disc 
  $\mathcal{D}_{\theta}(|\gamma|)$ defined by \eqref{D}
  for $\theta\in (-\pi, \pi)$.
  Since
  $1/|1+\frac{\xi}{\beta}| \leq 1$ for $\cos\, \theta \geq 0$ and
  $1/|1+\frac{\xi}{\beta}| \leq 1/|\sin \,\theta|$ for $\cos\, \theta < 0$ we find  that when $Re\,(\alpha) \leq 2$
   \ben
    \frac{1}{|(1+\frac{\xi}{\beta})^{k-\alpha+2}|}  \leq
     \left\{\begin{array}{ccc}
     	 A   & \textrm{for}  & \cos \,\theta \geq 0,\\[1.15ex]
     	 \frac{A}{|\sin \,\theta|^{k + 2 - Re\,(\alpha)}}  & \textrm{for}  & \cos \,\theta < 0
     \end{array}\right.
   \een   
   where $A=e^{-(Im \,\alpha)\,\arg (1+\frac{\xi}{\beta})} > 0$. 
   Thus in this case each integral can be analytically continued along any ray $\theta\neq \pi$ from $0$ to $+\infty\,e^{i \theta}$ and 
   defines a holomorphic function in the open disc
   $Re\, \left(\frac{e^{i \theta}}{x}\right) > 0$ whose opening is $< \pi$.  
   Let $x\in \mathcal{D}_{\theta}(|\gamma|)$ and let $\kappa=\arg (x)$. Note that from $x\in\mathcal{D}_{\theta}(|\gamma|)$ it follows
   that $|x| < \frac{\cos (\theta-\kappa)}{|\gamma|}$. Then we find that
   when $\cos \,\theta \geq 0$ and $Re (\alpha) \leq 2$ 
   $$\,
   \left|\sum_{k=0}^{\infty} \frac{\gamma^k\,x^k}{k!}
   \left(\int_0^{+\infty\,e^{i \theta}} 
   \frac{e^{-\frac{\xi}{x}}}{(1+\frac{\xi}{\beta})^{2-\alpha+k}}\,d \xi\right)\right| 
   \leq \frac{A}{c} 
   \sum_{k=0}^{\infty} \frac{|\gamma|^k\,|x|^k}{k!} <
      \frac{A}{c} 
      \sum_{k=0}^{\infty} \frac{\cos^k (\theta-\kappa)}{k!} < \infty\,,
   \,$$
   where $c=Re \left(\frac{e^{i \theta}}{x}\right) > |\gamma| > 0$. 
   Similarly, when $cos \,\theta < 0$  and $Re (\alpha) \leq 2$ we have that
   \ben
   \left|\sum_{k=0}^{\infty} \frac{\gamma^k\,x^k}{k!}
   \left(\int_0^{+\infty\,e^{i \theta}} 
   \frac{e^{-\frac{\xi}{x}}}{(1 + \frac{\xi}{\beta})^{2-\alpha+k}}\,d \xi\right)\right| 
   &\leq&  
   \frac{A}{c\,|\sin \theta|^{2- Re (\alpha)}} 
   \sum_{k=0}^{\infty} \frac{|\gamma|^k\,|x|^k}{k!\,|\sin\, \theta|^k} \\[0.2ex]
   &<&
   \frac{A}{c\,|\sin\ \theta|^{2-Re (\alpha)}} 
   \sum_{k=0}^{\infty} \frac{\cos ^k (\theta-\kappa)}{k!\,|\sin \theta|^k } < \infty\,.
   \een
   Thus from the Wiierstrass's theorem it follows that when $Re (\alpha) \leq 2$ the functional series
    \be\label{ser}
     \sum_{k=0}^{\infty} \frac{\gamma^k\,x^k}{k!}
     \left(\int_0^{+\infty\,e^{i \theta}} 
     \frac{e^{-\frac{\xi}{x}}}{(1+\frac{\xi}{\beta})^{2-\alpha+k}}\,d \xi\right)=
     \sum_{k=0}^{\infty} \frac{\gamma^k\,x^k}{k!}\,\varphi^{\theta}_k(x)
   \ee
    converges uniformly on the compact sets of the open disc $\mathcal{D}_{\theta}(|\gamma|)$.
    Since for all $k \geq 0$ the functions $\frac{\gamma^k\,x^k}{k!}\,\varphi^{\theta}_k(x)$ are  holomorphic functions
    on $\mathcal{D}_{\theta}(|\gamma|)$ so the sum.

   Assume now that there exists $k_1\in\NN_0$ such that 
   $k+2 - Re\,(\alpha) < 0$ for $0 \leq k \leq k_1$ while $k_1+1 +2 -Re (\alpha) \geq 0$. 
     Since $|(1+\frac{\xi}{\beta})^{\alpha-k-2}| \leq A\,(1+\frac{|\xi|}{\beta})^{Re (\alpha) -k-2}$
   where $A$ is above we find that for $0 \leq k \leq k_1$ 
    \ben
         & &
     \left|\int_0^{+\infty e^{i \theta}} 
       \frac{e^{-\frac{\xi}{x}}}{\left(1+\frac{\xi}{\beta}\right)^{2-\alpha+k}}\,d \xi\right| \leq
       A\,2^{Re (\alpha) -k-2} 
       \int_0^{\beta e^{i \theta}}
        e^{-|\xi|\,c}\,d \xi \\[0.2ex]
          &+&
          A\,\left(\frac{2}{\beta}\right)^{Re (\alpha) -k-2} 
          \int_{\beta\,e^{i \theta}}^{+\infty e^{i \theta}}
          |\xi|^{Re (\alpha) -k-2}\,e^{-c\,|\xi|}\,d \xi \\[0.2ex]
           &\leq&
         \frac{2^{Re (\alpha) -k-2}\,A}{c} (1-e^{-\beta\,c}) 
         	 +\frac{A}{c^{Re (\alpha)-k-1}}\,
         	\left(\frac{2}{\beta}\right)^{Re (\alpha) -k-2} 
         \,\Gamma(Re (\alpha)-k-1)\,.           
    \een
   Therefore for $0 \leq k \leq k_1$  we also can continue  analytically each integral along any ray $\theta\neq \pi$  and 
    define  holomorphic functions $\varphi^{\theta}_k(x)$ in the open disc
    $Re\, \left(\frac{e^{i \theta}}{x}\right) > 0$ whose opening is $< \pi$. Then for $x\in\mathcal{D}_{\theta}(|\gamma|)$
    we have  
    \ben
            & &
    \left|\sum_{k=0}^{\infty} \frac{\gamma^k\,x^k}{k!}
    \left(\int_0^{+\infty\,e^{i \theta}} 
    \frac{e^{-\frac{\xi}{x}}}{(1+\frac{\xi}{\beta})^{2-\alpha+k}}\,d \xi\right)\right| 
    \leq
   K \sum_{k=0}^{k_1} \frac{|\gamma|^k\,|x|^k}{2^k\,k!} \\[0.25ex]
            &+&
    \frac{A}{c^{Re (\alpha)-1}}\left(\frac{2}{\beta}\right)^{Re (\alpha)-2}
    \sum_{k=0}^{k_1} \frac{|\gamma|^k\,c^k\,\beta^k\,\Gamma(Re (\alpha) - k-1)}{2^k\,k!}\,x^k+ |F_{k_1+1}(x)|\\[0.25ex]
               &<&
        K \sum_{k=0}^{k_1} \frac{\cos^k (\theta-\kappa)}{2^k\,k!} +     
         \frac{A}{c^{Re (\alpha)-1}}\left(\frac{2}{\beta}\right)^{Re (\alpha)-2}
         \sum_{k=0}^{k_1} \frac{\cos^k (\theta-\kappa)\,c^k\,\beta^k\,\Gamma(Re (\alpha) - k-1)}{2^k\,k!}\\[0.25ex] 
           &+&
           |F_{k_1+1}(x)|    
    \een
    where $K=\frac{2^{Re (\alpha) -2}\,A}{c} (1-e^{-\beta\,c})$. 
    For the last addend we have that
   $$\,|F_{k_1+1}(x)| \leq \frac{A}{c} 
   \sum_{k_1+1}^{\infty} \frac{|\gamma|^k\,|x|^k}{k!} < 
   \frac{A}{c}    \sum_{k_1+1}^{\infty} \frac{\cos^k (\theta-\kappa)}{k!} < \infty 
   \,$$
    when $\cos \,\theta \geq$ and
    $$\,
    |F_{k_1+1}(x)| \leq \frac{A}{c\,|\sin\, \theta|^{2- Re (\alpha)}} 
    \sum_{k_1+1}^{\infty} \frac{|\gamma|^k\,|x|^k}{k!\,|\sin \,\theta|^k} <
    \frac{A}
    {c |\sin \theta|^{2- Re (\alpha)}}
       \sum_{k_1+1}^{\infty} \frac{\cos^k (\theta-\kappa)}{k!\,|\sin \,\theta|^k} < \infty
       \,$$
     when $\cos \,\theta < 0$.
   
      As a result  the functional series \eqref{ser}
  defines a holomorphic function on the disc $\mathcal{D}_{\theta}(|\gamma|)$ whose opening is $< \pi$.
        We denote this function by $\varphi_{\theta}(x)$.
        When we move $\theta\in (-\pi, \pi)$  the holomorphic functions $\varphi_{\theta}(x)$ glue together 
        analytically and define a holomorphic function $\tilde{\varphi}(x)$ on an open sector 
        $\widetilde{\mathcal{D}}$ from \eqref{D} whose  opening $> \pi$. 
        The restriction of $\tilde{\varphi}(x)$ on $\CC^*$ is a multivalued function.
        For every direction $\theta\neq \pi$ it has only one value. Near the singular direction $\theta=\pi$ it has two different values:
         $\varphi^{+}_{\theta}(x)=\sum_{k=0}^{\infty} \frac{\gamma^k\,x^k}{k!}\,
         \varphi^{\theta+\epsilon}_k(x)$ and
          $\varphi^{-}_{\theta}(x)=\sum_{k=0}^{\infty} \frac{\gamma^k\,x^k}{k!}\,
          \varphi^{\theta-\epsilon}_k(x)$ where $\epsilon > 0$ is a small number. 
                  
        This ends the proof.
    \qed

  \bre{dis}
   It seems that the function $\varphi_{\theta}(x)$ from \eqref{S} plays the part of 1 sum of the formal series
    $$\,
     \hat{\varphi}(x)=\sum_{k=0}^{\infty} 
     \frac{\gamma^k\,x^k}{k!}\,\hat{\varphi}_k(x)\,,
    \,$$ 
    where
     $$\,
      \hat{\varphi}_k(x)=\sum_{n=0}^{\infty}
       \frac{(-1)^n\,(2-\alpha+k)_n}{\beta^n}\,x^{n+1}\,.
     \,$$
     Here $(a)_n$ is the Pochhammer symbol.
     
     We are going to discuss in our next work the summation of such a series $\hat{\varphi}(x)$.
     
  \ere

   \bth{st}
    Assume that $\beta > 0$. Then with respect to the actual fundamental matrix solution $\Phi_0(x, 0)$
    defined by \thref{0-1} the DCHE \eqref{eq1} has a Stokes matrix $St_{\pi}$ at the origin in the form
     $$\,
      St_{\pi}=\left(\begin{array}{cc}
        1    & \mu\\
        0    & 1
         \end{array}\right)\,,
     \,$$
     where $\mu$ is introduced by \eqref{sm}. 
   \ethe

   \proof
   
   Let $\epsilon > 0$ be a small number and let $\theta=\pi$. Let $\Phi_0^{+}(x, 0)=\Phi_0^{\pi+\epsilon}(x, 0)$ and 
   $\Phi_0^{-}(x, 0)=\Phi_0^{\pi-\epsilon}(x, 0)$
   be the actual fundamental matrix solutions at the origin of the DCHE built by \thref{0-1}. To find the Stokes matrix $St_{\pi}$
    we have to compare the solutions $\Phi_0^{+}(x, 0)$   and $\Phi_0^{-}(x, 0)$
     $$\,
      \Phi_0^{-}(x, 0)=\Phi_0^{+}(x,0)\,
       \left(\begin{array}{cc}
          1    & \mu\\
          0    & 1
              \end{array}\right)\,.
     \,$$ 	 
   When $\alpha\in\NN$ we find that
     \ben
       \mu  &=&
       \frac{\gamma^{\alpha-1}\,e^{-\frac{\beta}{x}}}{\beta}\,
       	\left[\psi^{-}_{\pi}(x) - \psi^{+}_{\pi}(x)\right]=
       	2\,\pi\,i\,\gamma^{\alpha-1}\,e^{-\frac{\beta}{x}}\,
       	\res \left(\frac{v(\xi)\,e^{-\frac{\xi}{x}}}{\xi+\beta}\,;\,\xi=-\beta\right)\\[0.2ex]
       	 &=&
       	2\,\pi\,i\,\gamma^{\alpha-1}\,v(-\beta)=
       	2\,\pi\,i\,\gamma^{\alpha-1} 
       	\sum_{n=0}^{\infty} \frac{(-1)^n\,\beta^n\,\gamma^n}{n!\,\Gamma(\alpha+n)}\,. 
      \een
     Similarly, when $\alpha\in\ZZ_{\leq 0}$ we find that
       \ben
       \mu  &=&
       \frac{x^{1-\alpha}\,e^{-\frac{\beta}{x}}}{\beta}\,
       \left[\phi^{-}_{\pi}(x) - \phi^{+}_{\pi}(x)\right]=
       2\,\pi\,i\,e^{-\frac{\beta}{x}}\,
       \res \left(\frac{\xi^{1-\alpha}\,q(\xi)\,e^{-\frac{\xi}{x}}}{\xi+\beta}\,;\,\xi=-\beta\right)\\[0.2ex]
       &=&
       2\,\pi\,i\,(-\beta)^{1-\alpha}\,q(-\beta)=
       2\,\pi\,i\,(-\beta)^{1-\alpha} 
       \sum_{n=0}^{\infty} \frac{(-1)^n\,\beta^n\,\gamma^n}{n!\,\Gamma(2-\alpha+n)}\,. 
       \een
   Now we will show that the so found multipliers coincide. Indeed, let $\alpha\in\NN$. Then the multiplier
   $\mu$ corresponding to $\alpha\in\ZZ_{\leq 0}$ becomes
    $$\,
     \mu=2 \,\pi\,i\,(-\beta)^{1-\alpha}
     \sum_{n=\alpha-1}^{\infty}
     \frac{(-1)^n\,\beta^n\,\gamma^n}{n!\,(1+n-\alpha)}
    \,$$       
    since $1/\Gamma(z)=0$ for $z\in\ZZ_{\leq 0}$. Then
     $$\,
      \mu=2\,\pi\,i\,(-\beta)^{1-\alpha}
      \sum_{p=0}^{\infty}
      \frac{(-1)^{p+\alpha-1}\,\beta^{p+\alpha-1}\,\gamma^{p+\alpha-1}}
       {(p+\alpha-1)!\,p!}=
       2\,\pi\,i\,\gamma^{\alpha-1}
       \sum_{p=0}^{\infty}
       \frac{(-1)^p\,\beta^p\,\gamma^p}{p!\,\Gamma(p+\alpha)}\,,
     \,$$
     which is the multiplier $\mu$ corresponding to $\alpha\in\NN$.
     
     Let now $\alpha\notin \ZZ$ and let us compare the functions $\varphi^{\pi-\epsilon}_k(x)$ and $\varphi^{\pi+\epsilon}_k(x)$. We
     have
      $$\,
       \varphi^{\pi-\epsilon}_k(x) - \varphi^{\pi+\epsilon}_k(x)=
       (\beta)^{2-\alpha+k}\int_{\gamma} (\beta+\xi)^{\varepsilon-k-2}\,e^{-\frac{\xi}{x}}\,d \xi\,,
      \,$$
      where $\gamma=(\pi-\epsilon) - (\pi+\epsilon)$. Without changing the integral we can deform the path $\gamma$ into a Henkel
      type path going along the negative real axis from $-\infty$ to $-\beta$, encircling $-\beta$ in the positive sense and
      backing to $-\infty$. Then
       \ben
        \varphi^{\pi-\epsilon}_k(x) - \varphi^{\pi+\epsilon}_k(x) &=&
         \beta^{2-\alpha+k}\,\left(1-e^{-2 \pi\,i(\alpha-k-2)}\right)
         \int_{-\beta}^{-\infty} (\beta+\xi)^{\alpha-k-2}\,e^{-\frac{\xi}{x}}\,d \xi\\[0.2ex]
             &=&
          \beta^{2-\alpha+k}\, \left(1-e^{-2 \pi\,i\,\alpha}\right)\,e^{\frac{\beta}{x}}
           \int_0^{-\infty} u^{\alpha-k-2}\,e^{-\frac{u}{x}}\,d u\\[0.1ex]
            &=&
         \beta^{2-\alpha+k}\,\left(1-e^{-2 \pi\,i\,\alpha}\right)\,e^{\frac{\beta}{x}}\,x^{\alpha-k-1}\,
         \int_0^{+\infty} \tau^{\alpha-k-2}\,e^{-\tau}\,d \tau\\[0.15ex]
          &=&
           \beta^{2-\alpha+k}\,\left(1-e^{-2 \pi\,i\,\alpha}\right)\,\Gamma(\alpha-k-1)\,x^{\alpha-k-1}\,e^{\frac{\beta}{x}}\\[0.15ex]
             &=&
           -\frac{2\,\pi\,i\,(-1)^k\,e^{-\pi\,i\,\alpha}\,\beta^{2-\alpha+k}}{\Gamma(2-\alpha+k)}\,
           x^{\alpha-k-1}\,e^{\frac{\beta}{x}}\,,
           \een
        where we have used the Euler's reflection formula $\Gamma(1-z)\,\Gamma(z)=\frac{\pi}{\sin\,\pi z}$ for $z\notin\ZZ$.
        Then for the multiplier $\mu$ we find
         \ben
          \mu &=&
          \frac{x^{1-\alpha}\,e^{-\frac{\beta}{x}}}{\beta}\,
          \left[\varphi^{-}_{\pi}(x) - \varphi^{+}_{\pi}(x)\right]\\[0.2ex]
             &=&
         \frac{x^{1-\alpha}\,e^{-\frac{\beta}{x}}}{\beta}\,
         \sum_{k=0}^{\infty} \frac{\gamma^k\,x^k}{k!}\,
         \left[\varphi^{\pi-\epsilon}_k(x) - \varphi^{\pi+\epsilon}_k(x)\right]\\[0.2ex]
          &=&
          2 \pi\,i\,(-\beta)^{1-\alpha}
          \sum_{k=0}^{\infty} 
          \frac{(-1)^k\,\gamma^k\,\beta^k}{k!\,\Gamma(2-\alpha+k)}\,.
         \een
        
     This ends the proof.
	\qed

	 \vspace{1cm}
	 {\bf Acknowledgments.}	The author is grateful to L. Gavrilov for helpful discussions and comments.
	  The author was partially supported by Grant KP-06-N 62/5 \,
	  \textquotedblright{\it Algebraic and Analytic Methods in Differential Equations and Geometry}\textquotedblleft\, of the Bulgarian Found 
	  "Scientific research".

%%%%%%%%%%%%%%%%%%%%%%%%%%%%%%%%%%%%%%%%%%%%%%%%%%%%%%%
%apnedix
%%%%%%%%%%%%%%%%%%%%%%%%%%%%%%%%%%%%%%%%%%%%%%%%%%%%%%	
	\section*{Appendix}
	In this paragraph we will show by a direct computation that for lower values of $\alpha$ 
	 $$\,
	  \lim_{\sqrt{\varepsilon} \rightarrow 0} (q_R + q_L)=\lambda
	 \,$$
	where 
	\be\label{a-a}
	 \lambda=
	 \frac{\gamma^{\alpha-1}}{(\alpha-1)!}\,.
	 	\ee
	
	Let ${\bf \alpha=2}$. Then the multipliers $q_R$ and $q_L$ become
	\ben
	 q_R=\frac{1}{2 \sqrt{\varepsilon}}
	 \left(\frac{1+\varepsilon}{1-\varepsilon}\right)^{\frac{\gamma}{2 \sqrt{\varepsilon}}}\,
	 \quad
	  q_L=-\frac{1}{2 \sqrt{\varepsilon}}
	  \left(\frac{1-\varepsilon}{1+\varepsilon}\right)^{\frac{\gamma}{2 \sqrt{\varepsilon}}}\,.
	\een
	Note that in this case $ q_R \to +\infty$ and
       $ q_L \to -\infty$ when $\sqrt{\varepsilon} \to 0$. 
    Next the functions
    $\left(\frac{1-\varepsilon}{1+\varepsilon}\right)^{\frac{\gamma}{2 \sqrt{\varepsilon}}}$ and
    	 $\left(\frac{1+\varepsilon}{1-\varepsilon}\right)^{\frac{\gamma}{2 \sqrt{\varepsilon}}}$ are
   expressed as power series in $\sqrt{\varepsilon}$ as follows
   \be\label{sqrt}
   	 \left(\frac{1-\varepsilon}{1+\varepsilon}\right)^{\frac{\gamma}{2 \sqrt{\varepsilon}}} &=&
   	 \sum_{k=0}^{\infty}
   	 \left(-2 \gamma \sum_{p=o}^{\infty} \frac{(\sqrt{\varepsilon})^{4 p +1}}{4 p+2}\right)^k \frac{1}{k!}=
   	 1 - 2 \gamma \left(\frac{\sqrt{\varepsilon}}{2} + \frac{(\sqrt{\varepsilon})^5}{6} + \cdots\right)\nonumber\\[0.2ex]
   	   &+&
   	\frac{4 \gamma^2}{2 !}
   	\left(\frac{\sqrt{\varepsilon}}{2} + \frac{(\sqrt{\varepsilon})^5}{6} + \cdots\right)^2 -
   	\frac{8 \gamma^3}{3!}
   	\left(\frac{\sqrt{\varepsilon}}{2} + \frac{(\sqrt{\varepsilon})^5}{6} + \cdots\right)^3 + \cdots\\[0.2ex]
   		 \left(\frac{1+\varepsilon}{1-\varepsilon}\right)^{\frac{\gamma}{2 \sqrt{\varepsilon}}} &=&
   		 \sum_{k=0}^{\infty}
   		 \left(2 \gamma \sum_{p=o}^{\infty} \frac{(\sqrt{\varepsilon})^{4 p +1}}{4 p+2}\right)^k \frac{1}{k!}=
   		 1 + 2 \gamma \left(\frac{\sqrt{\varepsilon}}{2} + \frac{(\sqrt{\varepsilon})^5}{6} + \cdots\right)\nonumber\\[0.2ex]
   		 &+&
   		 \frac{4 \gamma^2}{2 !}
   		 \left(\frac{\sqrt{\varepsilon}}{2} + \frac{(\sqrt{\varepsilon})^5}{6} + \cdots\right)^2 +
   		 \frac{8 \gamma^3}{3!}
   		 \left(\frac{\sqrt{\varepsilon}}{2} + \frac{(\sqrt{\varepsilon})^5}{6} + \cdots\right)^3 + \cdots\nonumber
     \ee
     Then we find that when $\alpha=2$
      \ben
      \lim_{\sqrt{\varepsilon} \rightarrow 0} (q_R+q_L)=
      \lim_{\sqrt{\varepsilon} \rightarrow 0}
       \frac{2 \gamma\,\sqrt{\varepsilon} + O (\varepsilon)}{2 \sqrt{\varepsilon}}=\gamma\,,
      \een
      which coincides with \eqref{a-a}.

      Let ${\bf \alpha=4}$. Then the multipliers $q_R$ and $q_L$ become
      \ben
        q_R &=& -\left(\frac{1+\varepsilon}{1-\varepsilon}\right)^{\frac{\gamma}{2 \sqrt{\varepsilon}}}
        \frac{1} {(2 \sqrt{\varepsilon})^3}
        \left[\Gamma(3) - \Gamma(2)\,\frac{2 \sqrt{\varepsilon}\,\gamma}{1-\varepsilon^2}\right],\\[0.2ex]
            q_L &=& \left(\frac{1-\varepsilon}{1+\varepsilon}\right)^{\frac{\gamma}{2 \sqrt{\varepsilon}}}
              \frac{1} {(2 \sqrt{\varepsilon})^3} 
            \left[\Gamma(3) + \Gamma(2) \frac{ 2 \sqrt{\varepsilon}\,\gamma}{1-\varepsilon^2}\right]\,.
        	\een
        	Note that in this case $\lim_{\sqrt{\varepsilon} \rightarrow 0} q_R=-\infty$ and
         $\lim_{\sqrt{\varepsilon} \rightarrow 0} q_L=+\infty$. Again using the expressions \eqref{sqrt} we find
         that
         $$\,
          \lim_{\sqrt{\varepsilon} \rightarrow 0} (q_R+q_L)=
              \frac{4}{3}
          \lim_{\sqrt{\varepsilon} \rightarrow 0}
          \frac{(\sqrt{\varepsilon})^3 \,\gamma^3 + O (\varepsilon^2)}
          {(2 \sqrt{\varepsilon})^3}=
          \frac{\gamma^3}{3!}\,,
         \,$$	
         which coincides with \eqref{a-a}.

	  Let ${\bf \alpha=6}$. Then the multipliers $q_R$ and $q_L$ become
	  \ben
	   q_R &=&
	  \left(\frac{1+\varepsilon}{1-\varepsilon}\right)^{\frac{\gamma}{2 \sqrt{\varepsilon}}}
	   \frac{1}{2!\,(2 \sqrt{\varepsilon})^5}
	   \left[\frac{\Gamma(5)}{2!} -\Gamma(4)\,\frac{2\sqrt{\varepsilon}\,\gamma}{1-\varepsilon^2}
	   +\frac{\Gamma(3)}{2!}\,\left(\frac{2 \sqrt{\varepsilon}}{1-\varepsilon^2}\right)^2 \right],\\[0.2ex]
	     q_L &=&-
	     \left(\frac{1-\varepsilon}{1+\varepsilon}\right)^{\frac{\gamma}{2 \sqrt{\varepsilon}}}
	     \frac{1}{2!\,(2 \sqrt{\varepsilon})^5}
	     \left[\frac{\Gamma(5)}{2!} +\Gamma(4)\,\frac{2\sqrt{\varepsilon}\,\gamma}{1-\varepsilon^2}
	     +\frac{\Gamma(3)}{2!}\,\left(\frac{2 \sqrt{\varepsilon}}{1-\varepsilon^2}\right)^2 \right]\,.
	  \een
	  Then $\lim_{\sqrt{\varepsilon} \to 0} q_R=+\infty$ while $\lim_{\sqrt{\varepsilon} \to 0} q_L=-\infty$. 
	  For the sum $q_R + q_L$ we have
	   $$\,
	    \lim_{\sqrt{\varepsilon} \to 0} (q_R + q_L)=
	    \frac{4}{15} \lim_{\sqrt{\varepsilon} \to 0}
	    \frac{(\sqrt{\varepsilon})^5\,\gamma^5 + O (\varepsilon)^3}{(2 \sqrt{\varepsilon})^5}=\frac{\gamma^5}{5!}\,,
	   \,$$
	   which coincides with \eqref{a-a}.

	%%%%%%%%%%%%% References %%%%%%%%%%%%%%%%%%%%%%%%%%%%%%%%%%%%%%%%%%%%%%%
\begin{small}
    
\end{small}
%%%%%%%%%%%%%%%%%%%%%%%%%%%%%%%%%%%%%%%%%%%%%%%%%%%%%%%%%%%%%%%

\end{document}